\documentclass[a4paper,12pt]{article}
\usepackage{amsthm}
\usepackage{listings}
\usepackage{verbatim}
\usepackage{amsmath,amssymb}
\usepackage{amsfonts}
\usepackage{graphicx}
\usepackage{pdfpages}
\usepackage{caption}
\usepackage{subcaption}
\usepackage[english]{babel}
\usepackage{hyperref}
\usepackage{float}
\usepackage[T1]{fontenc}
\usepackage[utf8]{inputenc}
\usepackage{authblk}
\usepackage{multirow}
\usepackage{pdflscape}

\usepackage{perpage} 
\usepackage{chngcntr}
\theoremstyle{definition}
\usepackage[toc,page]{appendix}
\usepackage{tocloft}
\usepackage[all]{xy}
\usepackage{sectsty}
\usepackage[top=25mm, bottom=25mm, left=15mm, right=15mm]{geometry}

\begin{document}

\title{On the numerical accuracy of the method of multiple scales for nonlinear dispersive wave equations}
\author[*]{David Juhasz}
\author[*]{Per Kristen Jakobsen}
\affil[*]{Department of Mathematics and Statistics, the Arctic University of Norway, 9019 Troms\o, Norway}
\date{\today}

\maketitle
\begin{abstract}
In this paper we study dispersive wave equation using the method of multiple scales (MMS) and perform several numerical tests to investigate its accuracy. The key feature of our MMS solution is the linearity of the amplitude equation and the complex nature of the time-frequency $\omega$. The MMS is tested as an initial value problem using three choices of the dispersion model, one toy and two Lorentz models. Depending on the parameters of the problem, the amplitude equation can be both well- or ill-posed. Despite the ill-posedness, the MMS solution remains a valid approximation of the solution to the original nonlinear model.
\end{abstract}

\section{Introduction}

For mathematical models of waves propagating in material media, the phenomenon of \textit{dispersion} is frequently caused by the requirement of causality. 

For example, for the case of a light pulse propagating through an isotropic and homogeneous media, the propagating pulse induces a local dipole density, $\textbf{P}$, which for the simplest cases takes the form
\begin{align}
\textbf{P}(\textbf{x},t)&=\varepsilon_0\int_{-\infty}^t\mathrm{d}t'\chi(t-t')\textbf{E}(\textbf{x},t'),\label{eq0}
\end{align}
meaning that the polarization at a time $t$ only depends on the electric field at times previous
to $t$. This memory effect, which is the embodiment of causality,  is in optics called \textit{temporal dispersion}, or just dispersion. In this paper, in order to make our discussion specific, we will focus on this particular context in our work, but our methods and results apply quite widely to dispersive wave propagation.

The presence of dispersion evidently spells trouble for the integration
of the governing equations for the waves.  In general, they cannot be solved as an initial value problem.

This phenomenon is of course well known in the theory of wave propagation, and in optics in particular, and various more or less ingenious methods has been invented for getting around the problem. 


In optics one frequently tries to get around this problem by solving Maxwell's equations for optical pulse propagation as a boundary value problem, rather as an initial value problem.

In fact, one could argue that the boundary value problem is more closely aligned with the way experiments are done, than the initial value problem is. Waves in a material slab is launched by shining  a source laser at the interface of the slab, and therefore one could say that the incoming field at the boundary of the slab is fully controlled by using the laser. Thus, we have the data necessary for solving a boundary value  problem for optical pulse propagation. 

This is the basic, and in fact the only premise, which underlie the Unidirectional Pulse Propagation Equation(UPPE)~\cite{UPPE1}\cite{UPPE2}\cite{UPPE3} approach to optical pulse propagation. Of course, if there is significant back scatter of light from the interface and/or the material, the boundary condition is not fully controlled, and UPPE, and also other, less general, boundary solvers, which all rely on unidirectionality, are in trouble. For such cases, one can still solve Maxwell's equations for optical pulse propagation as a boundary value problem by using a more general approach than UPPE.  This approach is called the Bidirectional Pulse Propagation Equation(BPPE)~\cite{BPPE1}\cite{BPPE2}.

The BPPE approach is an exact method in the sense that no solutions has been lost when transitioning from Maxwell to BPPE. It is however also a purely numerical method and does not give any analytical insight into the pulse propagation problem. The UPPE is also an exact method, but only if one restricts to solutions of Maxwell than satisfy the condition of unidirectionality, and it is also a purely numerical method.

A much older approach to optical pulse propagation solves Maxwell's equations as an initial value problem by restricting to a class of solution that are spectrally narrow, also called narrow band solutions. This restriction makes it possible to derive equations, in general called {\it amplitude equations}, that, for a limited time, give a good approximation to Maxwell's equations, for solutions that are spectrally narrow. These amplitude equations can be solved as initial value problems. The systematic approach for deriving these equations is the method of multiple scales(MMS) \cite{Moloney}\cite{Per}. In optics, the best known such amplitude equation is the nonlinear Schr\o dinger equation(NLSE). Another well known amplitude equation, which is particularly useful for driven optical systems like a Laser, is the Complex Ginsburg-Landau equation~\cite{GinzLandau1}.

From a purely numerical point of view, the great thing about amplitude equations is that they are much faster to solve numerically than the original Maxwell's equations. The root cause of this is that for the narrow spectrum solutions represented by the amplitude equations, the fast frequency at the center of the spectrum, the so-called carrier wave,  needs not be temporally resolved, only deviations from the center frequency needs to be resolved, and the range of these deviations is by assumption small. 

In addition to being fast to solve numerically as an initial value problem, amplitude equations tend to have a universal form, at least to leading order, and this universal form, quite frequently, makes the equation amenable to analytical investigations. In the best cases, a complete analytical solution can be found. This is true for the NLSE equation \cite{Zakharov}. However, this analytic solvability is not robust. If we want to extend the amplitude equation beyond the leading order, which we must, if we want an equation that approximate the narrow spectrum solutions to Maxwell for a longer interval of time, the analytic solvability is typically lost, but the fast numerical solvability is not. 
So, one can say that the important feature of amplitude equations in optics,  is that they form the basis for a fast numerical approach for obtaining certain types of solutions to Maxwell's equations.

The aim of this paper is two-fold. 

Firstly, we want to go beyond the use of stationary modes as a basis for the MMS expansion, and consider the more general case of modes that are decaying, usually because of loss of some sort in the material, or gain for that matter. In this paper we will focus on the situation of lossy materials. The case of weak losses can treated by using stationary modes, assuming that the decaying terms are small compared to the leading part of Maxwell's equations. For situations where the loss is too large to be included as a perturbation, one must use a  MMS based on decaying modes. This is a situation that occur if one is investigating optical pulse propagation close to a material resonance. It is also the case if one want to derive amplitude equations in near-zero index situations. In fact, this last case is the major motivation for the work we do in this paper. Two interesting features of the resulting amplitude equations, features that we explore in detail later in this paper, is that the amplitude equations we derive are necessarily {\it linear } and that they are frequently {\it ill posed }. The interesting thing is that they, nevertheless, accurately represents both the linear and the nonlinear dynamics of the narrow band solutions of Maxwell's equations, for time intervals for which the amplitude equations according to the MMS procedure should approximate such solutions well. 

Secondly, we want to investigate the validity of the amplitude equations, as a fast numerical scheme for narrow band solutions to Maxwell's equations, by comparing the numerical solution of the amplitude equation to the corresponding numerical solution of Maxwell's equations. Since we cannot solve Maxwell's equation numerically as an initial value problem, this validity check has an obvious problem that needs to be handled. Handling this problem is the second major focus of this paper. It relies on the fact that in optics one almost always approximate the linear dispersion of materials using Sellmeier equations \cite{Sellmeier1}.




\numberwithin{equation}{section}
\section{A model wave equation, including nonlinearity and  general temporal dispersion}
The basic model equation we will use to illustrate our methods, and for which we will state our main conclusions, is the simplest nontrivial wave equation from nonlinear optics. It is scalar,  and includes general linear dispersion and a Kerr, cubic nonlinear material response. Our methods and conclusions apply much more widely than this, but for the sake of clarity and because numerical methods plays an important role in this paper, it is necessary to work within a specific class of equations.


In order to put our model equation into a real physical context, we start this section by deriving the equation from Maxwell's equations under some reasonable physical assumptions on the material response.

Maxwell's equations for a situation where there are no free charges or currents, are given by
\begin{align}
\partial_t\textbf{B}+\nabla\times\textbf{E}&=0,\nonumber\\
\partial_t\textbf{D}-\nabla\times\textbf{H}&=0,\nonumber\\
\nabla\cdot\textbf{D}&=0,\nonumber\\
\nabla\cdot\textbf{B}&=0.\label{eq1}
\end{align}
Most materials show no magnetic response at optical frequencies, thus we assume that
\begin{align}
\textbf{H}&=\frac{1}{\mu_0}\textbf{B},\nonumber\\
\textbf{D}&=\varepsilon_0\textbf{E}+\textbf{P}.\label{eq2}
\end{align}
The polarization is in general a sum of a term that is linear in $\textbf{E}$ and ones that are nonlinear in $\textbf{E}$. We thus have
\begin{align}
\textbf{P}=\textbf{P}_L+\textbf{P}_{NL},\label{eq3}
\end{align}
We assume that the linear material response is isotropic, homogeneous and causal
\begin{align}
\textbf{P}_L(\textbf{x},t)&=\varepsilon_0\int_{-\infty}^t\mathrm{d}t'\chi(t-t')\textbf{E}(\textbf{x},t').\label{eq4}
\end{align}

In the explicit calculations that we do in this paper, we will assume that the nonlinear polarization is restricted to the Kerr effect. Thus we will assume that
\begin{align}
\textbf{P}_{NL}=\varepsilon_0\eta\textbf{E}\cdot\textbf{E}\textbf{E},\label{eq5}
\end{align}
where $\eta$ is the Kerr coefficient. This is a choice we make just to be specific, the applicability of our methods, and the validity of our conclusions, derived in this paper, in no way depend on this particular choice for the nonlinear response.


Inserting (\ref{eq2})-(\ref{eq5}) into (\ref{eq1}), we can rewrite Maxwell's equations into the form
\begin{align}
\partial_t\textbf{B}+\nabla\times\textbf{E}&=0,\nonumber\\
\partial_t\textbf{E}-c^2\nabla\times\textbf{B}+\sqrt{2\pi}\partial_t\hat{\chi}(i\partial_t)\textbf{E}&=-\frac{1}{\varepsilon_0}\partial_t\textbf{P}_{NL},\nonumber\\
\nabla\cdot\left(\textbf{E}+\sqrt{2\pi}\hat{\chi}(i\partial_t)\textbf{E}\right)&=-\frac{1}{\varepsilon_0}\nabla\cdot\textbf{P}_{NL},\nonumber\\
\nabla\cdot\textbf{B}&=0,\label{eq7}
\end{align}
where we have used an alternative form of (\ref{eq4}) which is derived in Appendix A. The factor $\sqrt{2\pi}$ is a consequence of our conventions for the Fourier transform.

We will now restrict ourselves to solutions of the form
\begin{align}
\textbf{E}(z,t)&=E(z,t)\textbf{e}_y,\nonumber\\
\textbf{B}(z,t)&=B_1(z,t)\textbf{e}_x+B_2(z,t)\textbf{e}_z,\\
\textbf{P}_{NL}(z,t)&=P_{NL}(z,t)\textbf{e}_y,\label{eq8}
\end{align}
which are called {\it transverse electric waves}(TE). For this simplified case, Maxwell's equations take the form
\begin{align}
\partial_tB_1-\partial_zE&=0,\nonumber\\
\partial_tB_2&=0,\nonumber\\
\partial_tE-c^2\partial_zB_1+\sqrt{2\pi}\partial_t\hat{\chi}(i\partial_t)E&=-\frac{1}{\varepsilon_0}\partial_tP_{NL},\nonumber\\
\partial_zB_2&=0,\label{eq9}
\end{align}
where for the assumed Kerr effect we have
\begin{align}
P_{NL}&=\varepsilon_0\eta E^3.\label{eq10}
\end{align}
By taking cross derivatives it is easy to eliminate the magnetic field components and arrive at the equation
\begin{align}
\partial_{tt}E-c^2\partial_{zz}E+\sqrt{2\pi}\partial_{tt}\hat{\chi}(i\partial_t)E&=-\frac{1}{\varepsilon_0}\partial_{tt}P_{NL},\label{eq11}
\end{align}
which is the basic model equation we will be using in the rest of this paper.

\subsection{Scaling of the model equation}
If the aim is to solve Maxwell's equations numerically, in some specific physical context, it is not really necessary to scale the equation, and frequently this is not done, even if one can  argue that it could still be a useful thing to do. However, if one is going to derive an amplitude equation for the same physical situation, it might not be essential to scale the equation, but it certainly is extremely useful. After all, the essence of MMS is the ordering of terms in certain expansions according to size, and ensuring that this ordering, according to size, persists, up to some time of our choosing. 

Equation (\ref{eq11}) will be the staring point for our multiple scale approach. Let us start by picking scales $Z_0,T_0$ and $E_0$, for space, time and electric field amplitude, so that we have
\begin{align}
z&=Z_0 z',\nonumber\\
t&=T_0 t',\nonumber\\
E&=E_0 E',\label{eq12}
\end{align}
where the primed symbols are the scaled quantities. With these choices of scales the model equation (\ref{eq11}) takes the form


\begin{align}
\partial_{t't'}E'-\frac{c^2T_0^2}{Z_0^2}\partial_{z'z'}E'+\sqrt{2\pi}\partial_{t't'}\hat{\chi'}\left( i\partial_{t'}\right) E'&=-\eta E_0^2\partial_{t't'}E'^3,\label{eq19}
\end{align}
where $\hat{\chi'}\left( i\partial_{t'}\right)=\hat{\chi}\left( i\partial_{t'}\frac{1}{T_0}\right)$. We are at this point free to choose the time scale $T_0$, so let us choose it so that the factor before the $z$-derivative becomes one
\begin{align}
T_0=\frac{Z_0}{c}.\label{eq19.1}
\end{align}
 We next set the scale for the electric field to be the initial peak electric field amplitude. Thus
\begin{align}
    E_0=\underset{z}{max}\;|E(z,0)|.
\end{align}
With this we can write the model equation in the form

\begin{align}
\partial_{tt}E-\partial_{zz}E+\sqrt{2\pi}\partial_{tt}\hat{\chi}(i\partial_t)E&=-\varepsilon^2\partial_{tt}E^3,\label{eq20}
\end{align}
where we now have dropped primes on all quantities, since from this point on, only scaled quantities will appear. In this equation we have introduced the dimensionless parameter  $\varepsilon=\sqrt{\eta}E_0$. Typically, the Kerr parameter is fixed for any given material, whereas $E_0$ is at our disposal to vary over many order of magnitude, depending on the strength of the laser used to generate the initial electric field. This means that $\varepsilon$ can be made to vary over many orders of magnitude, but for realistic field intensities it is always smaller than one, usually much smaller than one.  $\varepsilon$ is the small perturbation parameter which we need for the MMS expansion.

In this paper, Fourier transforms and plane waves play a prominent role,  and whenever that is the case, it is convenient to pick the units for wave number and frequency in such a way that the phase of plane waves, and the Fourier transform, and its inverse, retain the same symbolic form in scaled and unscaled quantities. It is easy to verify that this is the case if we make the choice
\begin{align}
K_0&=\frac{1}{Z_0},\\
\Omega_0&=\frac{1}{T_0}.\label{eq20.1}
\end{align}
Thus we measure frequency in the well known unit Hertz, or cycles per unit time. Wavelength is in a similar way measured in periods per unit length. 

Note that (\ref{eq19.1}), which fixes the time scale in terms of the length scale, also, because of (\ref{eq20.1}), fix the frequency scale in terms of the wave number scale
\begin{align}
 \Omega_0=c K_0.\label{eq20.2}   
\end{align}
Thus, the only scale that remains to pick is the one for wave number. The initial field is in a lab situation generated using a laser. These days most optical labs have lasers that produce pulses of femto second duration~\cite{femtosecond1}, and labs with more  specialized equipment can produce pulses down to atto second durations~\cite{attosecond}. Such short pulses have a very broad spectrum and their dynamics are therefore hard to approximate using amplitude equations, which require narrow band pulses. There are versions of MMS that can handle such broad band pulses, but only in the weak dispersion limit. In this paper we apply MMS in a form that is tailored to the opposite limit of strong dispersion. Thus, in this paper we must assume that the initial field has a narrow wave number spectrum centred on a wave length determined by the lasing wave length of the laser generating the initial pulse. What the word "narrow" in the previous sentence means will be clarified in the MMS expansion in the next section. 

In this paper we are focused on validating amplitude equations derived using decaying modes in the vicinity of a material resonance, and it thus makes sense to pick the scale, $K_0$, for wave number, in such a way that the position of the resonance is centred on scales wavelengths that are of order one. The initial laser pulse will then have a  narrow band spectrum centred on a scaled wave length that is also of order one. We will in the rest of the paper assume that this has been done, and from now on only refer to scaled quantities, both in the model equation and in the specific material models that we will introduce in the sections to come.

\subsection{Decaying mode, amplitude equation, for the model equation}
We proceed with the multiple scale method by introducing the expansions
\begin{align}
\partial_t&=\partial_{t_0}+\varepsilon\partial_{t_1}+\varepsilon^2\partial_{t_2}+\ldots,\nonumber\\
\partial_z&=\partial_{z_0}+\varepsilon\partial_{z_1}+\varepsilon^2\partial_{z_2}+\ldots,\nonumber\\
e&=e_0+\varepsilon e_1+\varepsilon^2e_2+\ldots,\label{eq21}
\end{align}
where the connection between the multiple scale field amplitude $e$ and the electric field amplitude $E$ is given by
\begin{align}
E(z,t)&=e(z_0,t_0,z_1,t_1,\ldots)\big|_{t_j=\varepsilon^jt,z_j=\varepsilon^jz}.\label{eq21.1}  
\end{align}
The function $\hat{\chi}(i\partial_t)$ is expressed as a Taylor series as follows.
\begin{align}
\hat{\chi}(i\partial_t)&=\hat{\chi}(i(\partial_{t_0}+\varepsilon\partial_{t_1}+\varepsilon^2\partial_{t_2}+\ldots))=\hat{\chi}(i\partial_{t_0})+\hat{\chi}'(i\partial_{t_0})\left(\varepsilon i\partial_{t_1}+\varepsilon^2i\partial_{t_2}+\ldots\right)\nonumber\\
&+\frac{\hat{\chi}''(i\partial_{t_0})}{2}\left(\varepsilon i\partial_{t_1}+\varepsilon^2i\partial_{t_2}+\ldots\right)^2\nonumber\\
&=\hat{\chi}(i\partial_{t_0}) +\varepsilon i\partial_{t_1}\hat{\chi}'(i\partial_{t_0})+\varepsilon^2\left(i\partial_{t_2}\hat{\chi}'(i\partial_{t_0})-\frac{1}{2}\hat{\chi}''(i\partial_{t_0})\partial_{t_1t_1}+\ldots\right)+\ldots.\label{eq22}
\end{align}
We now insert (\ref{eq21}), (\ref{eq21.1}) and (\ref{eq22}) into (\ref{eq20}) and expand everything in sight. This gives us the following perturbation hierarchy
\begin{align}
\varepsilon^0:& &\partial_{t_0t_0}e_0-\partial_{z_0z_0}e_0+\sqrt{2\pi}\partial_{t_0t_0}\hat{\chi}(i\partial_{t_0})e_0&=0,\label{eq23}\\
\varepsilon^1:& &\partial_{t_0t_0}e_1-\partial_{z_0z_0}e_1+\sqrt{2\pi}\partial_{t_0t_0}\hat{\chi}(i\partial_{t_0})e_1&=-2\partial_{t_0t_1}e_0+2\partial_{z_0z_1}e_0-2\sqrt{2\pi}\partial_{t_0}\hat{\chi}(i\partial_{t_0})\partial_{t_1} e_0\nonumber\\
& & &-i \sqrt{2\pi}\partial_{t_0t_0}\hat{\chi}'(i\partial_{t_0})\partial_{t_1} e_0,\label{eq24}\\
\varepsilon^2:& &\partial_{t_0t_0}e_2-\partial_{z_0z_0}e_2+\sqrt{2\pi}\partial_{t_0t_0}\hat{\chi}(i\partial_{t_0})e_2&=-\partial_{t_1t_1}e_0-2\partial_{t_0t_2}e_0+\partial_{z_1z_1}e_0+2\partial_{z_0z_2}e_0\nonumber\\
& & &-\sqrt{2\pi}\hat{\chi}(i\partial_{t_0})\partial_{t_1t_1} e_0-2\sqrt{2\pi}\partial_{t_0}\hat{\chi}(i\partial_{t_0})\partial_{t_2} e_0\nonumber\\
& & &-2i \sqrt{2\pi}\partial_{t_0}\hat{\chi}'(i\partial_{t_0})\partial_{t_1t_1} e_0-i \sqrt{2\pi}\partial_{t_0t_0}\hat{\chi}'(i\partial_{t_0})\partial_{t_2} e_0\nonumber\\
& & &+\frac{1}{2}\sqrt{2\pi}\partial_{t_0t_0}\hat{\chi}''(i\partial_{t_0})\partial_{t_1t_1} e_0-\eta E_0^2\partial_{t_0t_0}e_0^3\nonumber\\
& & &-i \sqrt{2 \pi } \partial_{t_0t_0t_1} \hat{\chi}'(i\partial_{t_0})e_1-2 \sqrt{2 \pi } \partial_{t_0t_1} \hat{\chi}(i\partial_{t_0})e_1\nonumber\\
& & & -2 \partial_{t_0t_1}e_1+2 \partial_{z_0z_1} e_1.\label{eq25}
\end{align}
For the order $\varepsilon^0$ equation we choose the wave packet solution
\begin{align}
e_0(z_0,t_0,z_1,t_1,\ldots)=A_0(z_1,t_1,\ldots)e^{i\theta_0}+(*),\label{eq26}
\end{align}
where
\begin{align}
\theta_0=kz_0-\omega t_0,\label{eq27}
\end{align}
and where $\omega=\omega(k)$ is a {\it complex} solution to the dispersion relation
\begin{align}
\omega^2n^2(\omega)=k^2.\label{eq28}
\end{align}
Here, the complex refractive index, $n(\omega)$, is defined by
\begin{align}
n^2(\omega)=1+\sqrt{2\pi}\hat{\chi}(\omega).\label{eq29}
\end{align}
Observe that our multiple scale expansion is based on a complex, decaying mode, not a complex, stationary mode,  which is usual when one apply MMS far from any resonances of the material. We will see here and also in later sections that this fact will change many aspects of the resulting amplitude equations.

We must now calculate the right-hand side of the order $\varepsilon^1$ equation. Inserting (\ref{eq26}) into (\ref{eq25}) and (\ref{eq26}), we get
\begin{align}
\partial_{t_0t_0}e_2-\partial_{z_0z_0}e_2+\sqrt{2\pi}\partial_{t_0t_0}\hat{\chi}(i\partial_{t_0})e_2&=\left( 2i\omega\partial_{t_1}A_0+2ik\partial_{z_1}A_0+2i\omega\sqrt{2\pi}\hat{\chi}(\omega)\partial_{t_1} A_0\right.\nonumber\\
&\left.+i\omega^2\sqrt{2\pi}\hat{\chi}'(\omega)\partial_{t_1} A_0\right)e^{i\theta_0}+(*).\label{eq30}
\end{align}
In order to remove secular terms we must postulate that
\begin{align}
2ik\partial_{z_1}A_0+i\left( 2\omega+2\omega\sqrt{2\pi}\hat{\chi}(\omega)+\omega^2\sqrt{2\pi}\hat{\chi}'(\omega)\right)\partial_{t_1} A_0&=0.\label{eq31}
\end{align}
Observe that from the dispersion relation (\ref{eq28}) we have
\begin{align}
\omega^2(1+\sqrt{2\pi}\hat{\chi}(\omega))&=k^2,\nonumber\\
&\Downarrow\nonumber\\
\omega'(k)\left( 2\omega+2\omega\sqrt{2\pi}\hat{\chi}(\omega)+\sqrt{2\pi}\omega^2\hat{\chi}'(\omega)\right)&=2k.\label{eq32}
\end{align}
Thus (\ref{eq31}) can be written in the form
\begin{align}
\partial_{t_1}A_0+\omega'(k)\partial_{z_1}A_0=0,\label{eq33}
\end{align}
For the case of stationary modes this quantity is real and by definition equal to the group velocity for an initial light pulse with a narrow spectrum centred on the wave number $k$. 

The order $\varepsilon^1$ equation simplifies into
\begin{align}
\partial_{t_0t_0}e_1-\partial_{z_0z_0}e_1+\sqrt{2\pi}\partial_{t_0t_0}\hat{\chi}(i\partial_{t_0})e_1&=0.\label{eq34}
\end{align}
At this point we face a choice; which solution should we pick for this equation? The equation is homogeneous, and thus does not have any nontrivial particular solution, like the equation for $e_2$, at the next order, do. 

Here we pick the simplest possible solution 
\begin{align}
e_1=0,\label{eq35}
\end{align}
for (\ref{eq34}). The discussion of why we make this pick here, and what the consequences would be to make another less trivial choice, is best postponed until after we complete the derivation of the amplitude equation at order $\varepsilon^2$.

We now must compute the right-hand side of the order $\varepsilon^2$ equation. Inserting (\ref{eq35}) into the right-hand side of the order $\varepsilon^2$ equation we get
\begin{align}
\partial_{t_0t_0}e_2-\partial_{z_0z_0}e_2+\sqrt{2\pi}\partial_{t_0t_0}\hat{\chi}(i\partial_{t_0})e_2&=\left(-\partial_{t_1t_1}A_0+2i\omega\partial_{t_2}A_0+\partial_{z_1z_1}A_0+2ik\partial_{z_2}A_0\right.\nonumber\\
&-\sqrt{2\pi}\hat{\chi}(\omega)\partial_{t_1t_1}A_0+i2\omega\sqrt{2\pi}\hat{\chi}(\omega)\partial_{t_2}A_0\nonumber\\
&-2\omega \sqrt{2\pi}\hat{\chi}'(\omega)\partial_{t_1t_1}A_0+i\omega^2\sqrt{2\pi}\hat{\chi}'(\omega)\partial_{t_2}A_0\nonumber\\
&\left.-\frac{1}{2}\omega^2\sqrt{2\pi}\hat{\chi}''(\omega)\partial_{t_1t_1}A_0\right)e^{i\theta_0}+9\omega^2A_0^3e^{i3\theta_0}\nonumber\\
&+3(\omega_i-2\omega)^2|A_0|^2A_0e^{i\theta_0}e^{2t_0\omega_i}+(*),\label{eq36}
\end{align}
where $\omega_i=\text{Im} \;\omega$. At this point it is worth observing that none of the nonlinear terms in (\ref{eq36}) are secular. In addition to the usual nonsecular term $e^{i3\theta_0}$, we also have the term $e^{i\theta_0}e^{2t_0\omega_i}$, which would have been included into the secular terms, and thus, in the final amplitude equation, if it wasn't for the complex nature of $\omega$. This makes MMS based on decaying modes essentially different from the case of stationary modes.

Continuing the calculation, we observe that in order to remove secular terms from (\ref{eq36}), we must postulate that
\begin{align}
&-\partial_{t_1t_1}A_0+2i\omega\partial_{t_2}A_0+\partial_{z_1z_1}A_0+2ik\partial_{z_2}A_0-\sqrt{2\pi}\hat{\chi}(\omega)\partial_{t_1t_1}A_0+i2\omega\sqrt{2\pi}\hat{\chi}(\omega)\partial_{t_2}A_0\nonumber\\
&-2\omega \sqrt{2\pi}\hat{\chi}'(\omega)\partial_{t_1t_1}A_0+i\omega^2\sqrt{2\pi}\hat{\chi}'(\omega)\partial_{t_2}A_0-\frac{1}{2}\omega^2\sqrt{2\pi}\hat{\chi}''(\omega)\partial_{t_1t_1}A_0=0,\nonumber\\
&\Updownarrow\nonumber\\
&\partial_{t_1t_1}A_0\left(-1-\sqrt{2\pi}\hat{\chi}(\omega)-2\omega \sqrt{2\pi}\hat{\chi}'(\omega)-\frac{1}{2}\omega^2\sqrt{2\pi}\hat{\chi}''(\omega)\right)\nonumber\\
&+\partial_{t_2}A_0\left(2i\omega+i2\omega\sqrt{2\pi}\hat{\chi}(\omega)+i\omega^2\sqrt{2\pi}\hat{\chi}'(\omega)\right)+\partial_{z_1z_1}A_0+2ik\partial_{z_2}A_0=0.\label{eq37}
\end{align}
The factor multiplying the term $\partial_{t_2}A_0$ can be expressed using (\ref{eq32}), and equation (\ref{eq37}) therefore simplifies as follows
\begin{align}
&\partial_{t_1t_1}A_0\left(-1-\sqrt{2\pi}\hat{\chi}(\omega)-2\omega \sqrt{2\pi}\hat{\chi}'(\omega)-\frac{1}{2}\omega^2\sqrt{2\pi}\hat{\chi}''(\omega)\right)\nonumber\\
&+\partial_{t_2}A_0\frac{2ik}{\omega'(k)}+\partial_{z_1z_1}A_0+2ik\partial_{z_2}A_0=0,\nonumber\\
\nonumber\\
&\Downarrow\nonumber\\
\nonumber\\
&\partial_{t_2}A_0+\omega'(k)\partial_{z_2}A_0-i\beta\partial_{z_1z_1}A_0+i\alpha\partial_{t_1t_1}A_0=0,\label{eq38}
\end{align}
where
\begin{align}
\alpha&=\omega'(k)\frac{n^2(\omega)+2\omega \sqrt{2\pi}\hat{\chi}'(\omega)+\frac{1}{2}\omega^2\sqrt{2\pi}\hat{\chi}''(\omega)}{2k},\label{eq39}\\
\beta&=\frac{\omega'(k)}{2k}.\label{eq40}
\end{align}
By removing the secular terms from equation (\ref{eq36}), the order $\varepsilon^2$ equation turns into
\begin{align}
\partial_{t_0t_0}e_2-\partial_{z_0z_0}e_2+\sqrt{2\pi}\partial_{t_0t_0}\hat{\chi}(i\partial_{t_0})e_2&=9\omega^2A_0^3e^{i3\theta_0}\nonumber\\
&+3(\omega+2i\omega_i)^2|A_0|^2A_0e^{i\theta_0}e^{2t_0\omega_i}+(*).\label{eq41}
\end{align}
This equation is not homogeneous, and we chose at this point to solve for $e_2$, using only a particular solution. One such particular solution is evidently 
\begin{align}
e_2(z_0,t_0,\ldots)=c_1A_0^3e^{i3\theta_0}+c_2|A_0|^2A_0e^{i\theta_0}e^{2t_0\omega_i}+(*),\label{eq42}
\end{align}
where
\begin{align}
c_1&=\frac{\omega^2}{ k^2-\omega^2n^2(3\omega)}=\frac{1}{ n^2(\omega)-n^2(3\omega)},\label{eq43}\\
c_2&=\frac{3(\omega+2i\omega_i)^2}{k^2-\left(1+\hat{\chi}(\omega+i2\omega_i)\right) \left(\omega +2 i \omega_i\right) ^2}\label{eq44}
\end{align}
Defining an amplitude $A(z,t)$ by
\begin{align}
A(z,t)&=A_0(z_1,t_1,\ldots)\big|_{t_j=\varepsilon^jt,z_j=\varepsilon^jz},\label{eq45}
\end{align}
and proceeding in the usual way, using (\ref{eq33}) and (\ref{eq38}) we finally get the following amplitude equation 
\begin{align}
&\partial_{t}A+\omega'(k)\partial_{z}A-i\beta\partial_{zz}A+i\alpha\partial_{tt}A=0.\label{eq47}
\end{align}
The amplitude $A$ is related to the electric field amplitude $E$ through the formula
\begin{align}
  E(z,t)=A(z,t)e^{i(kz-\omega t)}+\varepsilon^2(c_1A^3(z,t)e^{i3(kz-\omega t)}+c_2|A(z,t)|^2A(z,t)e^{i(kz-\omega t)}e^{2t\omega_i})+(*).\label{eq47.1}
\end{align}
By design, for  (\ref{eq47}), (\ref{eq47.1}) to be an approximate solution to (\ref{eq20}), we must have
\begin{align}
\beta\partial_{zz}A\sim\alpha\partial_{tt}A\sim\mathcal{O}(\varepsilon^2),\nonumber\\
\partial_{t}A\sim\partial_{z}A\sim\mathcal{O}(\varepsilon),\label{eq48}
\end{align}
where we recall that $\varepsilon$ is a number much smaller than 1. Given these circumstances, we observe that
\begin{align}
\partial_{t}A&=-\omega'(k)\partial_{z}A\sim\mathcal{O}(\varepsilon),\nonumber\\
&\Downarrow\nonumber\\
\partial_{tt}A&=\omega'(k)^2\partial_{zz}A\sim\mathcal{O}(\varepsilon^2),\label{eq49}
\end{align}
and thus, the, second order in time, amplitude equation (\ref{eq47}), is asymptotically equivalent to the more convenient, first order in time, equation
\begin{align}
\partial_{t}A+\omega'(k)\partial_{z}A-i\left(\beta-\alpha\left(\omega'(k)\right)^2\right)\partial_{zz}A=0.\label{eq50}
\end{align}
This amplitude equation, together with relation (\ref{eq47.1}), are the two key elements defining a fast numerical scheme for narrow band solutions to (\ref{eq20}). 

Note that this amplitude equation is a linear equation. This is very different from the nonlinear Schr\"odinger equation, which is the leading order amplitude equation for our model equation far from any material resonances. The solution to the amplitude equation (\ref{eq50}), can be an accurate approximation our nonlinear model equation, despite the linearity of the amplitude equation, because the expression (\ref{eq47.1}) that connects the amplitude $A$ to the electric field amplitude $E$, is nonlinear. Note that if we want to have an amplitude equation which is a good approximation to the model equation beyond a time of order $\frac{1}{\varepsilon^2}$, we have to extend the MMS procedure to higher order in $\varepsilon$. It is evident from what we have said about secular terms in the paragraph following equation (\ref{eq1}), that these extended amplitude equations will all be linear, no matter to which order the MMS procedure is extended.

Since the amplitude equation is linear, the solution space, $S_A$, is of course a linear space. The relation (\ref{eq47.1}) amounts to a map, $M$, from $S_A$ into the solution space, $S_E$, of the exact equation  (\ref{eq20}). The map is certainly not surjective, and neither is it injective. The lack of injectivity means that the map cannot be used to induce a nonlinear superposition principle on its image set, $M(S_E)\subset S_E$, using the usual pullback/pushforward approach.

In order to get the approximate solution to our model equation defined by the amplitude equation (\ref{eq50}), and the relation (\ref{eq47.1}), we decided to pick particular solutions at order $\varepsilon$ and $\varepsilon^2$, at both orders disregarding the general solution to the homogeneous equation. The consequence of adding a solution to the homogeneous equation, in the form of wave packets,  at one or both orders, would be to add one or two new independent amplitudes to the problem. Each of these amplitudes would, through the removal of secular terms at order $\varepsilon$,  $\varepsilon^2$ and  $\varepsilon^3$, satisfy their own amplitude equations. Both these extra amplitude equations would also be linear, and uncoupled from each other and the one for the amplitude $A$. The relation defining the electric field in terms of the amplitudes would now be much more complicated and involve sums of products of all three amplitudes.

The deciding factor for whether we include these extra amplitudes, or not, is what kind of solutions of the model equation we are trying to approximate. This comes down to which kind of initial conditions for the model equation we are able to represent faithfully using our amplitude equation.

Any choice of a narrow-band initial amplitude for the amplitude equation (\ref{eq50}), will lead to a narrow band wave packet solution of the model equation, with dispersive properties determined by the choice of a complex solution $\omega=\omega(k)$ to the dispersion relation (\ref{eq28}), at order $\varepsilon^0$ of our MMS expansion. This is certainly a valid choice of initial condition, if our aim is to validate the amplitude equation (\ref{eq50}) together with its corresponding defining relation (\ref{eq47.1}) for $E$.

However, such a choice of initial condition for $E$ is awkward from a physical point of view. Usually  the initial field is determined by a laser, whose output, in most cases, can be approximated spectrally by a narrow Gaussian centred at the operating wave length of the laser. The initial condition for $E$, defined above using a narrow-band initial condition for $A$, consists of two Gaussians centred at $k$ and $3k$, and  with a special relationship between the amplitudes and phases of the two Gaussians. 

From the relation (\ref{eq47.1}), it is evident that in order for $E$ to be a Gaussian centred on some wave number $k$, the amplitude $A$ must be a Gaussian centred on $k=0$, and we have to introduce a solution to the homogeneous solution at order $\varepsilon^2$ of the form
\begin{align}
e_2(z_0,t_0,z_1,t_1,\ldots)=B_0(z_1,t_1,\ldots)e^{i(3kz-\omega(3k)t)}+(*),\label{eq51.1}
\end{align}
leading, in the end, to a relation determining the electric field in terms of the amplitudes $A$ and 
\begin{align*}
B(z,t)&=B_0(z_1,t_1,\ldots)\big|_{t_j=\varepsilon^jt,z_j=\varepsilon^jz},
\end{align*}
of the form
\begin{align}
  E(z,t)=A(z,t)e^{i(kz-\omega t)}+\varepsilon^2((B+c_1A^3(z,t))e^{i3(kz-\omega t)}\nonumber\\
  +c_2|A(z,t)|^2A(z,t)e^{i(kz-\omega t)}e^{2t\omega_i})+(*).\label{eq50.2}
\end{align}
By fixing the initial value of the amplitude $B$ to be 
\begin{align}
B(z,0)&=-c_1A^3(z,0)),\label{eq50.3}    
\end{align}
we ensure that the initial spectrum for $E$ is a Gaussian, whose center and width, is determined by the initial spectrum for the amplitude A.

The initial value for the amplitude $A$ is found by solving 
\begin{align}
  E(z,0)=A(z,0)e^{ikz}+\varepsilon^2
  c_2|A(z,0)|^2A(z,0)e^{ikz}+(*),\label{eq50.2.1}
\end{align}
for $A$ using for example Newton's method.

Since our aim is to get an approximation to the solution of the model equation to order $\varepsilon^2$ for times $t\le \varepsilon^{-2}$, we need to remove the secular terms generated by the amplitude $B$, at order $\varepsilon^3$,  since such terms would grow linearly and potentially disturb the spectrum at $3k$, already for times of order $\varepsilon^{-1}$.

From detailed calculations done while deriving the amplitude equation (\ref{eq50}), it is not hard to see that the resulting amplitude equation for the amplitude $B$ must take the form
\begin{align}
\partial_{t}B+\omega'(3k)\partial_{z}B=0.\label{eq50.4}
\end{align}
For the validity tests we will discussing later in the paper, we will always introduce such an extra  amplitude $B$, whenever it is necessary for  faithfully representing a Gaussian initial spectrum for $E$.


\section{Testing the validity of the amplitude equations}
The first goal of this paper was to use the MMS approach to derive the amplitude equation close to a material resonance for nonlinear wave equations with arbitrary, but of course causal, material response. Our main model equation and motivation comes from the field of nonlinear optics, but our methods and results can evidently be applied much wider than this.

Our second goal, which is the focus in this section and the rest of the paper, is to test the validity of the derived amplitude equations using high precision numerical simulations. If one is far from a material resonance, the dispersion relation of the material can to a good approximation be modelled by a polynomial function of the frequency, usually a second order polynomial is sufficient like in the Lorentz model~\cite{Lorentz}. In this situation, a direct simulation of the model equation can fairly easily be achieved, and the validity of the amplitude equation can thus be tested. However, even in this situation, numerical validation is not all that common. One reason for this is that the narrow band solutions, which are well approximated by the amplitude equations, are also the type of solutions to Maxwell that are most challenging to simulate numerically. This is because such solutions, which basically are wave packets, have a slowly varying, and consequently very wide, envelope, and at the same time contain a very large number of oscillation under the envelope. Thus one need both a large computational domain and a very high resolution of that domain. This makes for a large number discretization points and thus long running times.

When one is close to a material resonance, low order polynomial approximations does not work as well as when one is far from a resonance. In this situation different and more complicated approximations must be used. This typically turns the model equation, which is a differential equation in time, into an equation that is a  pseudo-differential equation in time. Solving such a thing as an initial value problem is not an easy proposition.

A frequently used class of such, more general approximations, are the rational functions. This type of approximations is in particular very much used for approximating the electric susceptibility in optics. For this situation, one can exactly transform the model equation into a equation that is a differential equation in time, and whose initial value problem can be solved numerically. The transformation is based on what is called the Sellmeier formulas in optics, and it thus make sense to denote the associated transformation for the {\it Sellmeier transformation}. 

\numberwithin{equation}{subsection}
\subsection{The Sellmeier Transformation}
The Sellmeier formulas in optics are approximations of the electric susceptibility in terms of sums rational functions of the simple type.
\begin{align}
  R(\omega)=\frac{1}{a\omega^2+b\omega+c},  \label{eq50.0}
\end{align}
where $a,b$ and $c$ are complex constants. Causality, and questions of loss or gain in the material, put restrictions on the constants that we leave aside for now. Later, when we do our numerical simulation in order to validate the amplitude equations, these restrictions come to the fore.
Finite sums of functions of the type (\ref{eq50.0}) will produce general rational functions. Thus the Sellmeier formulas are simply rational approximations to the electric susceptibility
\begin{align}
  \hat{\chi}(\omega)=\frac{P(\omega)}{Q(\omega)}.\label{eq50.01}    
\end{align}


In order to describe the Sellmeier transformation, observe first that  (\ref{eq20}) can be written in the form $\mathcal{L}(E,E^3)=0$, where $\mathcal{L}$ is a suitable operator that produces our equation (\ref{eq20}) and which includes the integral operator $\hat{\chi}(i\partial_t)$. Secondly, observe that the Fourier transform of our model equation can be rewritten into the form
\begin{align}
 \hat{\mathcal{L}}\left(E,E^3\right)&=0,\nonumber\\   
 \Updownarrow\nonumber\\
 \frac{1}{Q(\omega)}\left[Q(\omega)\hat{\mathcal{L}}\left(E,E^3\right)\right]&=0,\nonumber\\   
 \Updownarrow\nonumber\\
\frac{1}{Q(\omega)}\;\hat{\tilde{\mathcal{L}}}\left( E,E^3\right)&=0,\label{eq50.1}
\end{align}
where $\hat{\tilde{\mathcal{L}}}=Q(\omega)\hat{\mathcal{L}}$. 

Assuming the form of the susceptibility to be  described by the very simplest Sellmeier formula (\ref{eq50.0}), $\hat{\chi}(\omega)=R(\omega)$, the operator $\tilde{\mathcal{L}}$ is the differential operator 
\begin{align}
  \tilde{\mathcal{L}}&=\left(c+ib\partial_t-a\partial_t^2\right) \circ\mathcal{L}.\label{SellmeierOperator} 
\end{align}
The equation $\tilde{\mathcal{L}}\left(E,E^3\right)=0$ is the Sellmeier transformation of the original equation  $\mathcal{L}(E,E^3)=0$.

It is evident from (\ref{eq50.1}) that any solution to the differential equation 
\begin{align}
\tilde{\mathcal{L}}\left(E,E^3\right)&=0,\label{SellmeierTransformOfL}
\end{align}
is also a solution to the original pseudo differential equation (\ref{eq20}). If $S_M$ denotes the space of solutions for this original model equation and $S_D$ denote the space of solutions for the Sellmeier transformed equation (\ref{SellmeierTransformOfL}), we thus have $S_D\subset S_M$. We are not going to do a detailed mathematical analysis of situations where the Sellmeier transformation breaks down, these are situations where $Q(\omega)$ pass through zero, or is close to zero on a region of positive measure. This situation is very rarely  realized for ordinary materials.


Our idea is now to restrict the numerical validation to solutions in the smaller solution space $S_D$. This of course only make sense if the MMS procedure applied to the Sellmeier transformed equation (\ref{SellmeierTransformOfL}) produce the exact same amplitude equation as the one we got earlier, starting from the original pseudo differential equation (\ref{eq20}). From our short discussion of the relation between $S_M$ and $S_D$ in the previous paragraph, we certainly expect to get the same amplitude equation, but we still feel that it is prudent to directly verify that this is the case for each of the two explicit examples of Sellmeier formulas discussed in the following sections.

\numberwithin{equation}{subsection}
\subsection{A toy model for dispersion}
In this first validation calculation for our amplitude equation (\ref{eq50}) and associated reconstructed electric field amplitude (\ref{eq47.1}) we chose a material response function of the simple form
\begin{align}
\chi(t)=\left\{
\begin{array}{cc}
u e^{-v t}, & t>0\\
0, & t<0
\end{array}\right.,\label{eq51}
\end{align}
for some real positive constants $u$ and $v$. The corresponding electric susceptibility, which is the Fourier transform of (\ref{eq51}), is given by the formula 
\begin{align}
\hat{\chi}(\omega)&=\frac{u}{\sqrt{2\pi}}\frac{1}{v-i\omega}.\label{eq52}
\end{align}
Defining parameters $\gamma=\frac{v}{u},a=\frac{1}{u}$,  the formula for the susceptibility can be written in the more convenient form
\begin{align}
\hat{\chi}(\omega)&=\frac{1}{\sqrt{2\pi}}\frac{1}{\gamma-ia\omega}.\label{eq52.1}
\end{align}
We are not aware of any material response that in the optical regime is described by this electric susceptibility, but it does satisfy the, all important Kramer-Kronig relations, and thus describe a causal response. For us, this choice is a device for testing the amplitude equation in the simplest setting possible. Later in the paper we will investigate the validity of the amplitude equation when the electric susceptibility describe a very common atomic model, the Lorentz oscillator.  For this case the validation is conceptually the same as for the simple model above, but it is technically much harder.

\subsubsection{The Sellmeier transformation}
The model for electric susceptibility described above, leads to a Sellmeier transformed equation that, in our scheme of things, is as simple as possible to handle numerically. Following the procedure described in the previous section, we find that the Sellmeier transform of our model equation is 

\begin{align}
a\partial_{ttt}E+(\gamma+1)\partial_{tt}E-\gamma\partial_{zz}E-a\partial_{zzt}E+\varepsilon^2\gamma\partial_{tt}E^3+\varepsilon^2 a\partial_{ttt}E^3=0.\label{eq55}
\end{align}
This equation is, as expected, a differential equation and not a pseudo differential equation. It can be solved numerically using, for example a pseudo spectral method, which is what we do in this paper. As noted earlier, for the validation test to make sense, we must ensure that the MMS procedure applied to (\ref{eq55}) gives us the same amplitude equation (\ref{eq50}), which we got from the original model equation (\ref{eq20}). We do this is Appendix B and observe there that the derived amplitude equation is indeed the same as the original one (\ref{eq50}).

\subsubsection{Numerical results}
Actually solving  the Sellmeier transformed equation (\ref{eq55}), is awkward,  because it is not solved explicitly with respect to the highest derivative. However, the fact that $\varepsilon$ is small means that, by iteration, we can approximate (\ref{eq55}) by an equation that {\it is} solved explicitly with respect to highest derivative. We achieve this by expanding out the offending term $a\partial_{ttt}E^3$ as follows 
\begin{align}
a\partial_{ttt}E^3&=a\left(6\left(\partial_{t}E\right)^3+18E\partial_tE\partial_{tt}E+3E^2\partial_{ttt}E\right),\label{eq97}
\end{align}
and substituting for $\partial_{ttt}E$ using (\ref{eq55}).  

Dropping terms of order $\varepsilon^4$, which we must do in order to be consistent with our order $\varepsilon^2$  MMS expansion,  we get the following  explicit equation
\begin{align}
&a\partial_{ttt}E+(\gamma+1)\partial_{tt}E-\gamma\partial_{zz}E-a\partial_{zzt}E+\eta E_0^2\gamma\left(6E\left(\partial_{t}E\right)^2+3E^2\partial_{tt}E\right)\nonumber\\
&+\varepsilon^2a\left[6\left(\partial_{t}E\right)^3+18E\partial_tE\partial_{tt}E+\frac{3E^2}{a}\left(\gamma \partial_{zz}E+a\partial_{zzt}E-(\gamma+1)\partial_{tt}E\right)\right]=0.\label{eq98}
\end{align}

Since the explicit equation (\ref{eq98}) agree with the, exact, implicit equation (\ref{eq55}) to order $\varepsilon^2$, the amplitude equation to order $\varepsilon^2$ for these two equations must be the same. Thus we can use (\ref{eq98}) to validate our original amplitude equation (\ref{eq50}).

Recall that the amplitude equation is derived using one specific mode of the linearization. Thus, to do the validation we must select one of these modes. Given this mode, all the parameters in the amplitude equation are determined in terms of the parameters of (\ref{eq98}). The detailed expressions are derived in Appendix B.
This linear mode, which is an explicit function of $z$ and $t$, can now be used to determine all the three initial conditions needed to solve (\ref{eq98}) numerically, and thus to complete the validation.




The equation (\ref{eq98}) has three independent modes. The solution to the linear part of (\ref{eq98}) can be expressed as inverse Fourier integrals of all three modes 
\begin{align}
E(z,t)=\frac{1}{\sqrt{2\pi}}\left(\int_{-\infty}^\infty\mathrm{d}k A_1(k)e^{i(kz-\omega_1(k)t)}+\int_{-\infty}^\infty\mathrm{d}k A_2(k)e^{i(kz-\omega_2(k)t)}+\int_{-\infty}^\infty\mathrm{d}k A_3(k)e^{i(kz-\omega_3(k)t)}\right).\label{eq88}
\end{align}
The electric field $E(z,t)$ must be real so (\ref{eq88}) must be equal to its complex conjugate
\begin{align}
E^*(z,t)
&=\frac{1}{\sqrt{2\pi}}\left(\int_{-\infty}^\infty\mathrm{d}k A_1^*(-k)e^{i(kz+\omega_1^*(-k)t)}\right.\nonumber\\
&\left.+\int_{-\infty}^\infty\mathrm{d}k A_2^*(-k)e^{i(kz+\omega_2^*(-k)t)}+\int_{-\infty}^\infty\mathrm{d}k A_3^*(-k)e^{i(kz+\omega_3^*(-k))t)}\right).\label{eq89}
\end{align}
Let us first assume that $k>0$. From the form of the dispersion relation (\ref{eq62}) 
\begin{align}
p(k,\omega)&=a\omega^3+i\omega(\gamma+1)-ak^2\omega-i\gamma k^2,\label{eq89.1}
\end{align}
we observe that if $\omega(k)$ is a solution to $p(k,\omega)=0$, then we have 
\begin{align}
p(k,-\omega^*)&=-a(\omega^*)^3-i\omega^*(\gamma+1)+ak^2\omega^*-i\gamma k^2\nonumber\\
&=-\left[ a\omega^3+i\omega(\gamma+1)-ak^2\omega-i\gamma k^2\right]^* =0.\label{eq89.2}
\end{align}
This implies that if $\omega(k)$ is a solution, then $-\omega^*(k)$ is a solution as well. Thus, assuming that the three solutions are distinct, we can number them in such a way that 
\begin{align}
\omega_1^*(k)&=-\omega_2(k),\label{eq90}\\
\omega_3^*(k)&=-\omega_3(k).\label{eq91}
\end{align}
Next observe that the dispersion relation (\ref{eq62}) is even in $k$. Using this fact, we can number the solutions for negative $k$ is such a way that $\omega_{1,2,3}(k)=\omega_{1,2,3}(-k)$. Using these relations, together with (\ref{eq90}) and (\ref{eq91}), formula (\ref{eq89}) turns to
\begin{align}
E^*(z,t)&=\frac{1}{\sqrt{2\pi}}\left(\int_{-\infty}^\infty\mathrm{d}k A_1^*(-k)e^{i(kz-\omega_2(k)t)}+\int_{-\infty}^\infty\mathrm{d}k A_2^*(-k)e^{i(kz-\omega_1(k)t)}\right.\nonumber\\
&\left.+\int_{-\infty}^\infty\mathrm{d}k A_3^*(-k)e^{i(kz-\omega_3(k))t)}\right).\label{eq92}
\end{align}
For (\ref{eq88}) and (\ref{eq92}) to be the same, the following relations between the amplitudes $A_{1,2,3}(k)$ must hold
\begin{align}
A_1(k)&=A_2^*(-k),\nonumber\\
A_3(k)&=A_3^*(-k).\label{eq93}
\end{align}
\noindent
From these relations, we see that the amplitudes for $k<0$ are determined from their values for $k>0$ and vice versa.

While deriving the asymptotic solution in Appendix B, we assumed the solution to the 0-th order equation in the perturbation hierarchy (\ref{eq57}) to be one of the three possible modes. The amplitude equation (\ref{eq50}) is the correct equation only for narrow band solutions consisting of one such mode, centred around wave number $k_0$ which we, without loss of generality, can assume is positive. To be specific, assume that this mode is $A_1(k)$. The narrow band property implies that $A_1(k)=0$ for $k<0$, and from this the relations (\ref{eq93}) implies that we can consistently choose $A_3(k)=0$ for all $k$ and $A_2(k)=0$ for $k>0$. 

Using these assumptions and (\ref{eq93}), we get from (\ref{eq92})
\begin{align}
E(z,t)&=\frac{1}{\sqrt{2\pi}}\left(\int_{-\infty}^\infty\mathrm{d}k A_1^*(-k)e^{i(kz-\omega_2(k)t)}+\int_{-\infty}^\infty\mathrm{d}k A_2^*(-k)e^{i(kz-\omega_1(k)t)}\right)\nonumber\\
&=\frac{1}{\sqrt{2\pi}}\int_{-\infty}^\infty\mathrm{d}k A_1(k)e^{i(kz-\omega_1(k)t)}+\int_{-\infty}^\infty\mathrm{d}k A_1^*(k)e^{-i(kz-\omega_1^*(k)t)}\nonumber\\
&=\frac{1}{\sqrt{2\pi}}\int_{-\infty}^\infty\mathrm{d}k \underbrace{\left[A_1(k)e^{-i\omega(k)t}+A_1^*(-k)e^{i\omega^*(-k)t}\right]}_{\hat{E}(k,t)}e^{ikz}=\frac{1}{\sqrt{2\pi}}\int_{-\infty}^\infty\mathrm{d}k\hat{E}(k,t)e^{ikz},\label{eq95}
\end{align}
where we let $\omega_1(k)=\omega(k)$ to match the parameters in our asymptotic solution (\ref{eq47.1}) with (\ref{eq95}).

We now turn back to equation (\ref{eq98}). Upon transforming this PDE  into its spectral domain $k$, we get a third order ODE's with a nonlinear right-hand side. We want to solve this ODE  as an initial value problem, and therefore need three initial conditions, here denoted by $f(k),g(k)$ and $h(k)$. These three initial conditions we obtain from (\ref{eq95}) 
\begin{align}
\hat{f}(k)&=\hat{E}(k,0),\nonumber\\
\hat{g}(k)&=\partial_{t}\hat{E}(k,0),\nonumber\\
\hat{h}(k)&=\partial_{tt}\hat{E}(k,0).\label{eq95.1}
\end{align}
These conditions can be also expressed in terms of the amplitude $A_1(k)$ using (\ref{eq95})
\begin{align}
\hat{f}(k)&=A_1(k)+A_1^*(-k),\label{eq100}\\
\hat{g}(k)&=-i\omega(k)A_1(k)+i\omega^*(-k)A_1^*(-k),\label{eq101}\\
\hat{h}(k)&=-\omega^2(k)A_1(k)-\left(\omega^*(-k)\right)^2 A_1^*(-k),\label{eq102}
\end{align}
\noindent
The only thing left to do now, is to compute the initial condition for the amplitude equation (\ref{eq50}) in terms of our chosen amplitude $A_1(k)$. 

As discussed at the end of section 2.3, the initial condition for $E$, which is natural from a physical point of view, is one whose spectrum is a narrow Gaussian. This is taken care of by letting the spectral amplitude $A_1(k)$  be a Gaussian centered at some wave number $k_0$
\begin{align}
A_1(k)=\left\{
\begin{array}{cc}
De^{-\delta(k-k_0)^2},& k>0\\
0,& k<0 
\end{array}\right.,\label{eq104}
\end{align}
where $D,\delta>0$.

As also discussed at the end of section 2.3, we might have to introduce an extra amplitude at order $\varepsilon^2$ in order to faithfully represent the initial condition on $E$ in terms of an initial condition for the amplitude equation.

The way to determine if any extra amplitude has to be introduced at order $\varepsilon^2$, is to assume the opposite. Given this, the relation between the amplitude $A$ and the electric field $E$ is determined by (\ref{eq47.1}).

 We now take the inverse Fourier transform of equation (\ref{eq100}) and equate it to the right-hand side of (\ref{eq47.1}) evaluated at $t=0$. Matching separately the first part and the second part, which is the complex conjugates of the first part,  on both sides, we get
\begin{align}
\mathcal{F}^{-1}\left\{ A_1(k)\right\}&=A(z,0)e^{ik_0z}+c_1\eta E_0^2 A^3(z,0)e^{i3k_0z}+c_2\eta E_0^2|A(z,0)|^2A^*(z,0)e^{ik_0z},\nonumber\\
\mathcal{F}^{-1}\left\{ A_1^*(-k)\right\}&=A^*(z,0)e^{-ik_0z}+c_1^*\eta E_0^2\left( A^*\right)^3(z,0)e^{-i3k_0z}+c_2^*\eta E_0^2|A(z,0)|^2A(z,0)e^{-ik_0z},\label{eq103}
\end{align}
which is a nonlinear system of algebraic equations, consisting of two equations and two unknowns $A(z,0)$ and $A^*(z,0)$. This system can easily be solved numerically, for example using Newtons method.

Note that the parameter $\delta$ controls the width of $A_1(k)$. Given that the amplitude $A(z,t)$ should be slowly varying in $z$, $\partial_zA\sim \mathcal{O}(\varepsilon)$, the parameter $\delta$ should be chosen such that $\delta\sim 1/\varepsilon$.

\begin{figure}[t!]
  \centering
\captionsetup{width=0.85\textwidth}
\begin{subfigure}{.5\textwidth}
  \centering
  \includegraphics[scale=0.25]{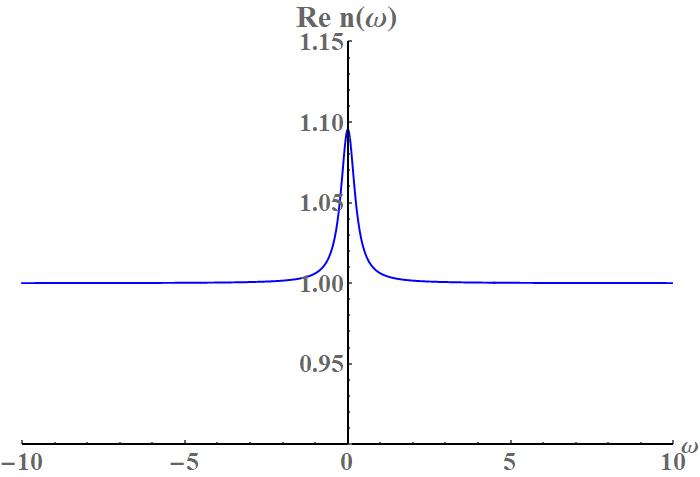}
  \caption{}
  \label{fig1a}
\end{subfigure}%
\begin{subfigure}{.5\textwidth}
  \centering
  \includegraphics[scale=0.25]{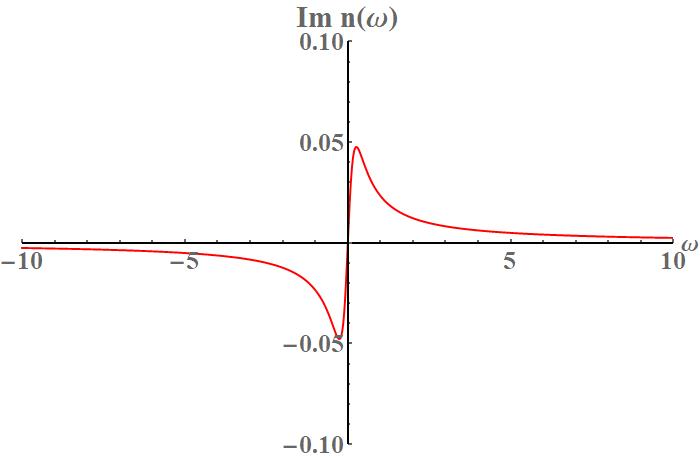}
  \caption{}
  \label{fig1b}
\end{subfigure}
\caption{Refractive index $n(\omega)$ for the toy model after scaling. The parameters in this figure are $\gamma=5$ and $a=20$.}
\label{fig1}
\end{figure}

In order to do a numerical comparison, we fix the refractive index by choosing the parameter values $\gamma=5$ and $a=20$. In figure \ref{fig1} we see the real and imaginary part of the resulting refractive index of our material.

At this point, the only remaining parameter to set is the, all important, nonlinearity parameter $\varepsilon$. In all our derivation, we assumed that this parameter was small, only then will there be separation of scales, which is the key assumption underlying our derivations of amplitude equations using MMS.

Usually, one is deriving amplitude equations in order to gain some analytical insight into the wave propagation problem at hand, and the actual numerical value of $\varepsilon$ does not need to be fixed. Here however, where we are doing a numerical validation of the amplitude equations, the situation is different, here we need to give a specific value to the nonlinearity parameter. And, with respect to this, we face a tradeoff. 

On the one hand, we should pick a value for $\varepsilon$ that is as small as possible, in order for the separation of scales to be as large as possible. Only for such values can we expect our amplitude equations to accurately approximate the original model equation.

On the other hand, a very small value for the nonlinearity parameter, means that we have to propagate the waves for very long distances for the nonlinearity in the model equation to influence the spectrum of the waves. Very long distances means very long times which translates into very long running time for the simulations of the model equation.

In this paper we choose the value for the nonlinearity parameter to be $\varepsilon=10^{-1}$. This certainly does not seem to be a very small number, but as we will see, even for a nonlinearity parameter as large as this, our amplitude equation does actually approximate the exact model equation well. This is just another example of what one could call the phenomenon of the {\it unreasonable accuracy of asymptotic methods}.

As it turns out, our chosen value for the nonlinearity parameter is also physically reasonable. For example, for a light pulse in a visible part of the spectrum, at $586\;nm$, propagating through Argon gas at atmospheric pressure, of an intensity equal to half the critical ionization threshold, the value for the nonlinearity parameter is  $\varepsilon=0.13$~\cite{argon}.

Let us now point our attention to calculating the correct initial condition for $A$. In order for the amplitude equation (\ref{eq50}) to be valid, the amplitude $A$ must be spectrally narrow. In fact, the whole MMS expansion is based on the assumption that the width of the spectrum of $A$ is of order epsilon. 

\begin{figure}[H]
  \centering
\captionsetup{width=0.85\textwidth}
  \includegraphics[scale=0.25]{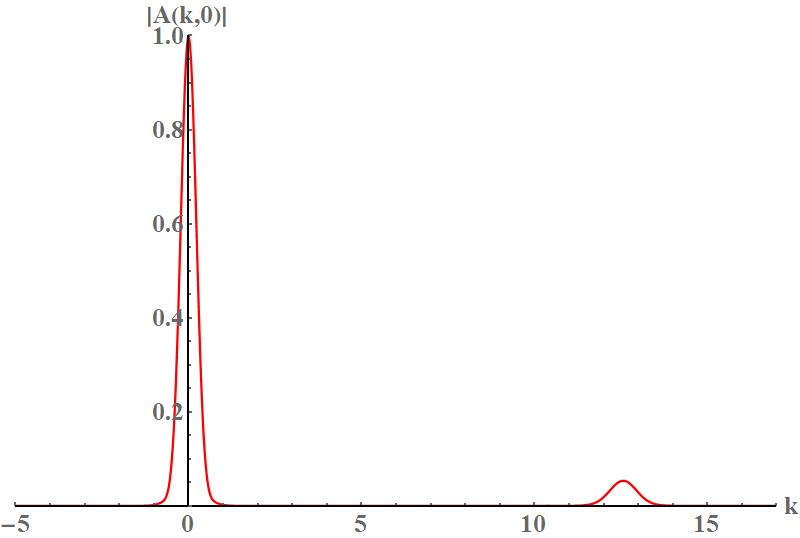}
  \caption{The Fourier transform of the initial conditions $A(k,0)$ obtained with Newtons iterative method from the equations (\ref{eq103}).}
\label{fig1.01}
\end{figure}

The spectrum for $A$ which we find assuming that no new amplitude is required at order $\varepsilon^2$ is displayed in figure \ref{fig1.01}. The parameters used to calculate the initial condition for $A$ were $a=20$ and $\gamma=5$.

Evidently, the width of this spectrum is not of order $\varepsilon$. The main part of the spectrum is narrow, but there is an additional peak in the spectrum centred at $2k_0$ whose distance from the origin is of order one, not order $\varepsilon$. Thus, assuming that no new amplitude is needed at order $\varepsilon^2$ leads to a contradiction.

From relation (\ref{eq47.1}) we observe that the peak in the spectrum for $A$ corresponds to a peak at $3k_0$ for $E$. Furthermore, from figure \ref{fig1.01} we observe that the height of the peak is of order $\varepsilon^2$. From these two observations it is evident that the peak in the spectrum for $A$ can be taken into account by introducing an extra amplitude $B$ in the MMS expansion at order $\varepsilon^2$. As we outlined at the end of section 2.3, the amplitude $B$ will come equipped with its own linear amplitude equation, decoupled from the one for $A$. 

The initial values for $A$ and $B$ are now found by solving equation (\ref{eq50.2.1}) for $A$, using the approach from (\ref{eq103}), and then using the identity (\ref{eq50.3}) to determine the initial value for $B$ in terms of the one for $A$.

\begin{figure}[H]
  \centering
\captionsetup{width=0.85\textwidth}
  \includegraphics[scale=0.3]{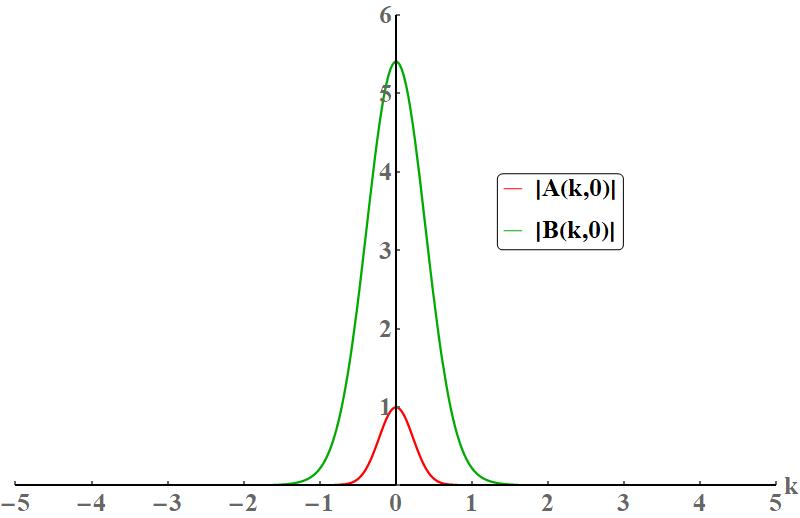}
  \caption{The Fourier transform of the initial conditions $A(k,0)$ and $B(k,0)$ obtained with Newtons iterative method from the equations (\ref{eq103}) and (\ref{eq50.3}).}
\label{fig1.1}
\end{figure}

In figure \ref{fig1.1} we display an example of the initial condition $A(k,0)$. In this example we use the value $k_0=2\pi$. Note that this means that our initial pulse spectrum is located to the right of the frequency defining the material resonance, in a region of anomalous dispersion,  which means that the real part of the refractive index decrease for increasing frequency. 



 Now we have everything we need in order to compare the numerical solution of the model equation (\ref{eq98}) to it's corresponding amplitude equation (\ref{eq50}).
\begin{figure}[H]
  \centering
\captionsetup{width=0.85\textwidth}
\begin{subfigure}{.5\textwidth}
  \centering
  \includegraphics[scale=0.33]{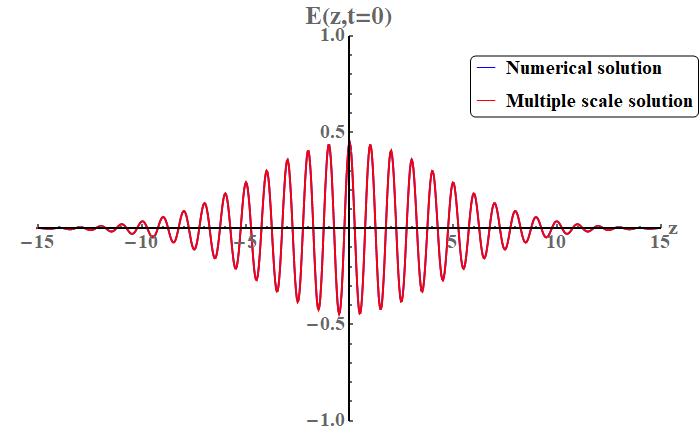}
  \caption{}
  \label{fig2a}
\end{subfigure}%
\begin{subfigure}{.5\textwidth}
  \centering
  \includegraphics[scale=0.33]{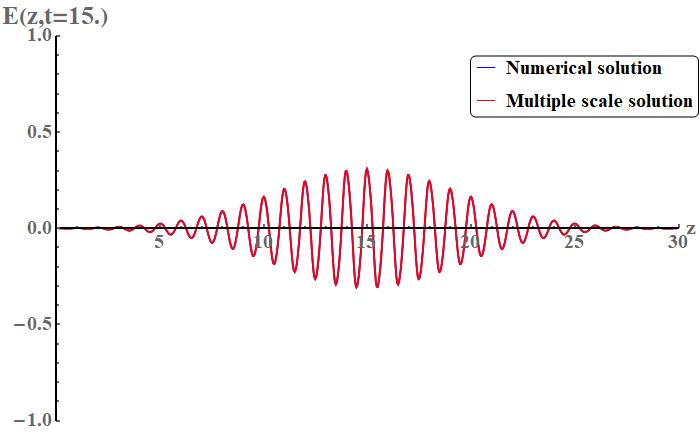}
  \caption{}
  \label{fig2b}
\end{subfigure}
\begin{subfigure}{.5\textwidth}
  \centering
  \includegraphics[scale=0.33]{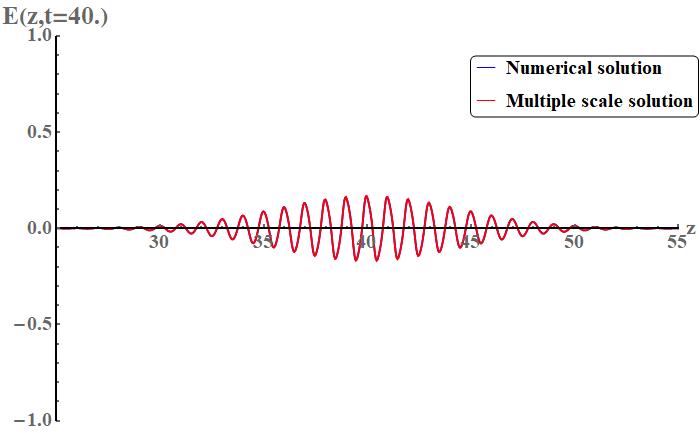}
  \caption{}
  \label{fig2c}
\end{subfigure}%
\begin{subfigure}{.5\textwidth}
  \centering
  \includegraphics[scale=0.33]{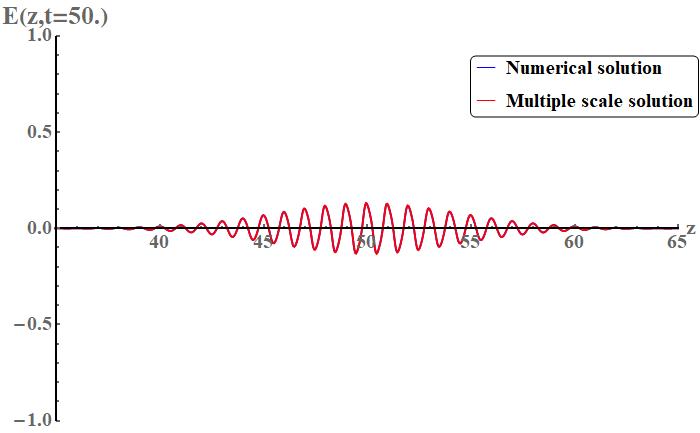}
  \caption{}
  \label{fig2d}
\end{subfigure}%
\caption{The solution of (\ref{eq98}) using a system of first order ODE's (blue graph) compared with solution (\ref{eq47.1}) using the amplitude $A(z,t)$ computed from (\ref{eq50}) (red graph) for the toy model dispersion.}
\label{fig2}
\end{figure}
\noindent
Figure \ref{fig2} depicts these two solutions at different times.  We can see that they are indeed very close. For the values of the parameters used in this comparison, the complex frequency of the chosen mode is given by $\omega_0=\omega(k_0)=6.28 - i2.5\times 10^{-2}$.

\begin{figure}[t!]
  \centering
\captionsetup{width=0.85\textwidth}
\begin{subfigure}{.5\textwidth}
  \centering
  \includegraphics[scale=0.25]{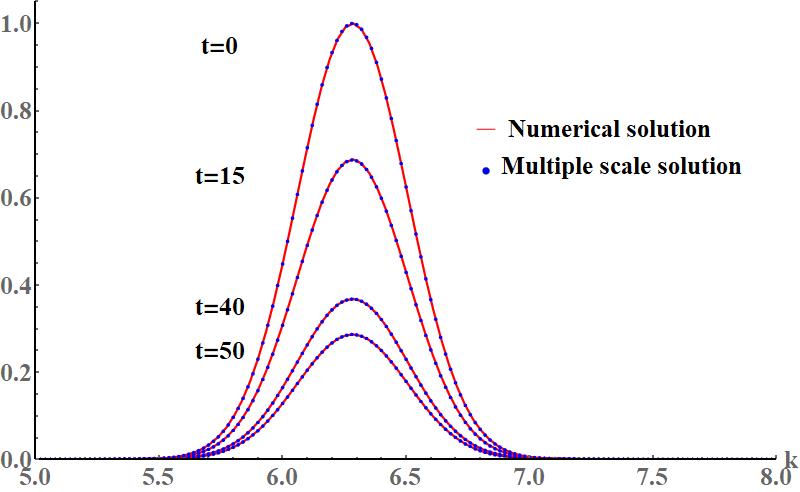}
  \caption{}
  \label{fig2.2a}
\end{subfigure}%
\begin{subfigure}{.5\textwidth}
  \centering
  \includegraphics[scale=0.25]{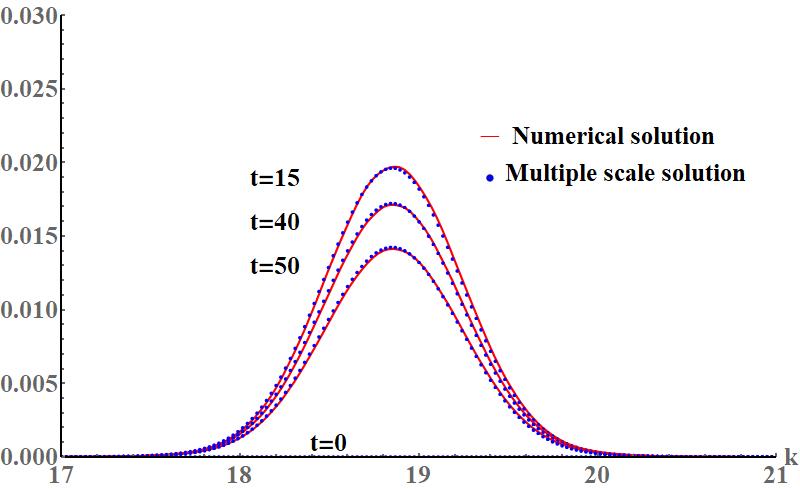}
  \caption{}
  \label{fig2.2b}
\end{subfigure}
\caption{The spectrum of the solutions of (\ref{eq98}) using a system of first order ODE's (blue graph) compared with the spectrum of the solution (\ref{eq47.1}) using the amplitude $A(z,t)$ computed from (\ref{eq50}) (red graph) for the toy model dispersion.}
\label{fig2.2}
\end{figure}

On order to get more insight into the accuracy, let  us compare the wave number spectrum of these two solutions. This is what is displayed in figure figure \ref{fig2.2}, which  includes both the major bumps in the $k$ spectrum. The larger one has its support at the frequency $k=k_0=2\pi$ as expected. This represents the linear part of the solution. The smaller one is located around the frequency $k=3k_0=6\pi$. It makes sense that it sits precisely at $3k_0$, since the Kerr nonlinearity has the form $E^3$, which means that when inserting a plane wave for $E$ into $E^3$, the wave numbers get multiplied by a factor of 3.

From figure \ref{fig2.2b} it does appear that there is a small deviation between the exact numerical solution and the solution derived from the amplitude equation.This deviation is however one or more orders of magnitude smaller than $\varepsilon^2$, and such deviations are to be expected for for a MMS expansion truncated at order $\varepsilon^2$.  

Based on the experience from testing the MMS solution, it became clear that in order to stay in the correct asymptotic regime, where the MMS solution is valid, the parameter values fixing the problem are subject to certain constraints.

First of all, the constants $\alpha$ and $\beta$ occurring in the amplitude equation, defined in (\ref{eq39}), (\ref{eq40}), include the derivatives of the susceptibility $\hat{\chi}(\omega)$. These derivatives come from the Taylor expansion of this function. We expect this Taylor series to converge, so by assumption, these derivatives must not break the order of the preceding terms in the expansion. This assumption depends mainly on the parameters $\gamma$ and $a$ in (\ref{eq52}) and the frequency $\omega$ around which the Taylor series is expanded. To preserve the order of these terms is closely tied to the possibility of making the assumptions in (\ref{eq48}) and consequently arriving at the amplitude equation. The assumption made in (\ref{eq48}) are related also to the choice of the parameter $\delta$ in the initial condition (\ref{eq104}). The amplitude $A(z,t)$ is by assumption slowly varying in $z$ compared with the exponential $e^{ik_0z}$ because of (\ref{eq45}). Therefore the initial condition (\ref{eq104}) should be fast varying in $k$ and the parameter $\delta$ needs to be chosen accordingly, for example $\delta=1/\varepsilon$.

\noindent

\subsubsection{Stability and well-posedness}
From the numerical solution of the amplitude equation in the previous section, and its ability to accurately approximate the exact, narrow band, solutions for Maxwell, one might think that our amplitude equation is just fine. However, nothing could be further from the truth. In this section we will show that the amplitude equation is in fact ill posed as a PDE.

\begin{figure}[t!]
  \centering
\captionsetup{width=0.85\textwidth}
  \includegraphics[scale=0.3]{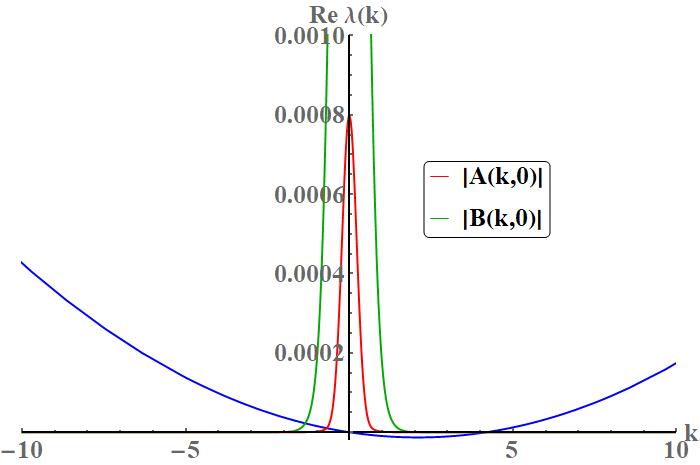}
  \caption{Stability of the amplitude equation (\ref{eq50}) for the toy model dispersion. The red graph is representing the amplitude $|A(k,0)|$. This function was scaled to fit the figure. The parameter values in this figure are $a_1=3\times 10^{-6}$ and $a_2=-1.26\times 10^{-5}$.}
\label{fig7}
\end{figure}

Since the amplitude equation is linear, this fact can be easily proven by merely calculating the growth curve for the equation. This curve we find by  inserting $A(z,t)=e^{\lambda(k) t}e^{ik z}$ into the equation, cancelling the common factor $e^{\lambda(k) t}e^{ik z}$, and extracting the real part of the resulting algebraic equation.
\begin{align}
\lambda+ik \omega'(k_0)-i\left(\beta-\alpha\left(\omega'(k_0)\right)^2\right)\left(ik\right)^2&=0,\nonumber\\
&\Updownarrow\nonumber\\
\lambda(k)&=-ik \omega'(k_0)-ik^2\left(\beta-\alpha\left(\omega'(k_0)\right)^2\right),\nonumber\\
&\Updownarrow\nonumber\\
\text{Re}\;\lambda(k)&=a_1k^2+a_2k,\label{eq105.2}
\end{align}
where $a_1$ and $a_2$ are real parameters that depend on the complex parameters $\alpha,\beta$ and $\omega'(k)$ is the derivative of the dispersion relation (\ref{eq28}). This function is a parabola which pass through the origin, and is displayed in figure \ref{fig7}, using the parameter values from our numerical test.

From figure \ref{fig7} it is evident that the amplitude equation is ill posed by definition; arbitrary high spatial frequencies will grow exponentially in time, with no upper bound for the growth rate.

Being ill posed is usually, for good reasons, thought of as ringing the death knell for any proposed mathematical model. But, still, in the situation investigated in the previous section, the ill posed amplitude equation is an excellent tool for simulating narrow band solutions to Maxwell's equation.

We regard this as another example of the way asymptotic methods makes use of, and gives meaning to, otherwise meaningless mathematical expressions. Anyone that is a user of asymptotic methods can not help noticing this fact. The classical case of this is Eulers example, where the useless power series $\sum_n(-1)^n n!\; x^n$, with zero radius of convergence, is, nevertheless, an accurate approximation to a certain exponential integral. Similarly, the expansion for the energy of the anharmonic quantum oscillator is an asymptotic series with zero radius of convergence~\cite{Anharmonic}. In fact, essentially all series used in quantum theory, in both the particle and the field incarnations of the theory, is known to, or believed to, have zero radius of convergence~\cite{ZeroConvergence1}\cite{ZeroConvergence2}. But, nevertheless, their usefulness is indisputable, their predictions give some of the most accurate correspondences between theory and experiment in all of science.

The reason why the ill posed amplitude equation nevertheless is an accurate numerical model for the narrow band solutions of Maxwell, for which it was designed, is simple. From the form of the stability curve (\ref{eq105.2}) it is clear that it passes through 0 for every parameter. Combining this with that fact that the initial condition for $A$ will always be centered at $k=0$, we can conclude that in the case of ill-posedness, the growth rate around $k=0$ will not be large. For any solution that satisfies the assumptions used to derive the amplitude equation, the fastest growing wave number component in the spectrum of the amplitude $A(z,t)$, will, during the time for which the amplitude equation is valid, $t\le \varepsilon^{-2}$, not grow large enough to affect the solution to order $\varepsilon^2$ or greater.

\subsection{Lorentzian model of dispersion}
In the previous section we used the simplest possible rational approximation to the electric susceptibility. As we noted there, this model is, as far as we know, not a realistic physical model, even though it is causal and thus does satisfy the optical  Kramer-Kronig relations. In this subsection we consider the next simplest one, which is of the form 
\begin{align}
\hat{\chi}(\omega)=\frac{1}{\sqrt{2\pi}}\frac{\omega_p^2}{\omega_r^2-\omega^2-i\gamma\omega}.\label{eq158}
\end{align}
This dispersion is a realistic physical model of the electric susceptibility. It can be derived from a purely classical model of an atom consisting of a positive charge, representing the nucleus, together with all the electrons, save one~\cite{Masud}. The remaining electron is singled out by being the one that resonantly respond to an imposed oscillatory electric field. Since the single electron is much lighter than the rest of the atom, which after all contains the nucleus, we in effect describing a simple oscillator, which is a harmonic oscillator, unless the field is strong enough to  pull the electron  too far from the atom.  The model is called the Lorentz oscillator, and have been applied to a vast range of materials in the gaseous, liquid and solid phase. The parameters $\omega_r,\omega_p$ and $\gamma$ are interpreted as the resonance frequency of the oscillator, the plasma frequency and the absorptive loss. The factor of $2\pi$ is, as noted before,  a consequence of our Fourier transform conventions, which are introduced in Appendix A. The Lorentz oscillator model is causal and thus satisfy the Kramer-Kronig relations.

We now repeat the calculations from the previous section for the Lorentz oscillator. We start by finding a  differential equation $\tilde{\mathcal{L}}\left(E,E^3\right)=0$, corresponding to our model pseudo differential equation (\ref{eq20}), using the Sellmeier transformation. After that, we use MMS to verify that we get the same amplitude equation for the Sellmeier transformed equation as the one we got from the original equation (\ref{eq20}). We then move on to doing a numerical comparison of the accuracy of the amplitude equation with respect to the Sellmeier transformed equation. As noted earlier in this paper, this amounts to a direct comparison of the numerical accuracy of the amplitude equation with respect to the original model equation  (\ref{eq20}), which is the goal of this paper.
\begin{figure}[t!]
  \centering
\captionsetup{width=0.85\textwidth}
  \includegraphics[scale=0.4]{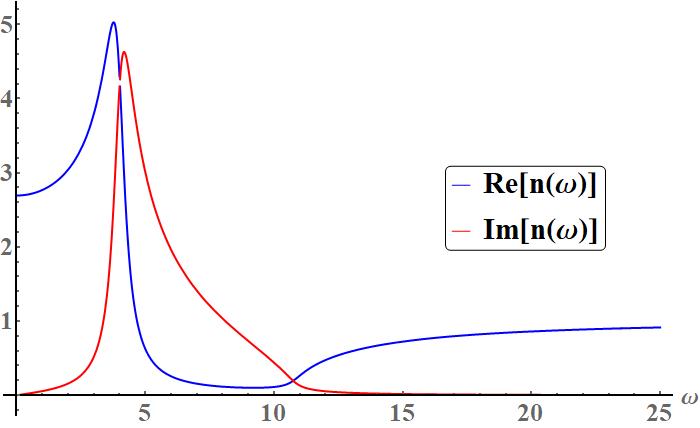}
  \caption{Refractive index $n(\omega)$ for the Lorentzian model of dispersion (\ref{eq158}) after scaling. The parameters in this figure are $a=-0.01,b=-7\times 10^{-3}$ and $c=0.16$.}
\label{fig3}
\end{figure}
In our calculations we rewrite the function (\ref{eq158}) in a for us more convenient form
\begin{align}
\hat{\chi}(\omega)=\frac{1}{\sqrt{2\pi}}\frac{1}{\frac{\omega_r^2}{\omega_p^2}-\frac{1}{\omega_p^2}\omega^2-i\frac{\gamma}{\omega_p^2}\omega}=\frac{1}{\sqrt{2\pi}}\frac{1}{a\omega^2+ib\omega+c},\label{eq160}
\end{align}
where we have defined $a=-1/\omega_p^2,
b=-\gamma/\omega_p^2$ and $c=\omega_r^2/\omega_p^2$. In our case, the refractive index is defined as $n^2(\omega)=1+\sqrt{2\pi}\hat{\chi}(\omega)$. In Figure (\ref{fig3}) we display an example of such a refractive index, corresponding to the numerical values $a=-0.01,b=-7\times 10^{-3}$ and $c=0.16$.

In the next subsections we derive the equation (\ref{eq50.1}) for the Lorentzian model of dispersion and obtain the MMS solution for it. Next we test the MMS solution using two choices of parameters for the Lorentz model.

\subsubsection{The Sellmeier transformation}
The Lorentzian model of dispersion leads to the following Sellmeier transformed equation
\begin{align}
&(1+c)\partial_{tt}E-b\partial_{ttt}E-a\partial_{tttt}E+(b\partial_{zzt}-c\partial_{zz}+a\partial_{zztt})E+\varepsilon^2(c\partial_{tt}-b\partial_{ttt}-a\partial_{tttt})E^3=0.\label{eq164}
\end{align}
Again, the derived equation (\ref{eq164}) does not include any pseudo differential operator. We proceed in this case in the same way as in chapter 3.2.1 and obtain MMS solutions to (\ref{eq164}) and see if they match with the solutions in the asymptotic regime for (\ref{eq22}). A detailed calculation of this MMS solution is done in Appendix C. The derived amplitude equation is found to be the same as (\ref{eq50}). 

\subsubsection{Numerical results}
We proceed by testing the results based on the amplitude equation (\ref{eq50}) and a solution to (\ref{eq164}) for the Lorentz model of dispersion. Two tests will be conducted, each of them with a different choice of parameters $a,b,c$ for the Lorentz model. 

The Sellmeier transformed equation (\ref{eq164}) is implisit in its highest time derivative, but can be approximated by an equation explicit in its highest time derivative. Following the same procedure as in chapter 3.2.2 we have
\begin{align}
&(1+c)\partial_{tt}E-b\partial_{ttt}E-a\partial_{tttt}E+(b\partial_{zzt}-c\partial_{zz}+a\partial_{zztt})E\nonumber\\
&=\varepsilon^2\left[-c\left(6E\left(\partial_{t}E\right)^2+3E^2\partial_{t}^{(2)}E\right)+b\left(6\left(\partial_{t}E\right)^3+18E\partial_tE\partial_{t}^{(2)}E+3E^2\partial_{t}^{(3)}E\right)\right.\nonumber\\
&\left.+a\left(3 \partial_{t}^{(4)}E E^2+18 E \left(\partial_{t}^{(2)}E\right)^2+24 \left(\partial_{t}^{(3)}E\right) E \partial_t E+36 \left(\partial_tE\right)^2 \partial_{t}^{(2)}E\right)\right],\label{eq218}
\end{align}
where $\partial_{t}^{(4)}E$ is expressed as
\begin{align}
\partial_t^{(4)}E&=\frac{1}{a}\left((1+c)\partial_{tt}E-b\partial_{ttt}E+(b\partial_{zzt}-c\partial_{zz}+a\partial_{zztt})E\right).\label{eq219}
\end{align}

As before, the amplitude (\ref{eq50}) equation is an equation for one of the four independent modes equation (\ref{eq218}) has. In order to test the MMS solution we need to choose one mode for which we have the amplitude equation (\ref{eq50}). In order to arrive at the point where we can chose the mode, we proceed like we did in subsection 3.2.2.
The electric field $E(z,t)$ is expressed as the inverse Fourier transform of the sum of all 4 modes and is then equated with its complex conjugate. In order to proceed we need to establish relations between the mode frequencies found by solving the dispersion equation  (\ref{eq170}). 

We observe that
\begin{align}
p(k,\omega)&=a \omega ^4+i b \omega ^3+\omega ^2(c+1)-k^2 \left(a \omega ^2+i b \omega +c\right)=0\nonumber\\
&\Downarrow\nonumber\\
p(k,-\omega^*)&=a(\omega^*) ^4-i b (\omega^*) ^3+(\omega^*) ^2(c+1)-k^2 \left(a (\omega^*) ^2-i b \omega^* +c\right)=0\nonumber\\
&=\left[ a \omega ^4+i b \omega ^3+\omega ^2(c+1)-k^2 \left(a \omega ^2+i b \omega +c\right)\right]^* =0\label{eq219.1}
\end{align}
Thus like in subsection 3.2.2, solution space of $p(k,\omega)=0$ is the same as for $p(k,-\omega^*)=0$. We can therefore conclude that if $\omega(k)$ is a solution, then $-\omega^*(k)$ is as well. Assuming that all solutions are distinct, we can enumerate the four the solutions in such a way that the following relations between $\omega_{1,2,3,4}(k)$ hold:
\begin{align}
\omega_1^*(-k)&=-\omega_2(k),\label{eq205}\\
\omega_3^*(-k)&=-\omega_4(k).\label{eq206}
\end{align}
The equation $p(k,\omega)=0$ is also even in $k$, and thus we can enumerate the four solutions corresponding to negative $k$ in such a way that $\omega_i(-k)=\omega_j(k)$, for some $i$ and $j$. Using this fact, together the relations (\ref{eq205}), (\ref{eq206}), the reality of the electric field implies the following relations between the amplitudes $A_{1,2,3,4}(k)$:
\begin{align}
A_1(k)&=A_2^*(-k),\label{eq209}\\
A_3(k)&=A_4^*(-k).\label{eq210}
\end{align}
The amplitudes for negative argument are defined from their values for $k>0$. As before, we want the electric field consisting of only one amplitude, for which we obtained the amplitude equation (\ref{eq50}), therefore we set $A_2(k)=A_3(k)=A_4(k)=0$ for $k>0$. This implies $A_{3,4}=0,\forall k$. We thus arrive at
\begin{align}
E(z,t)
&=\frac{1}{\sqrt{2\pi}}\int_{-\infty}^\infty\mathrm{d}k \underbrace{\left[A_1(k)e^{-i\omega(k)t}+A_1^*(-k)e^{i\omega^*(-k)t}\right]}_{\hat{E}(k,t)}e^{ikz}=\frac{1}{\sqrt{2\pi}}\int_{-\infty}^\infty\mathrm{d}k\hat{E}(k,t)e^{ikz}.\label{eq212}
\end{align}
The amplitude $A_1(k)$ can be chosen arbitrary and it will define the initial condition $E(z,0)$.

The equation (\ref{eq218}) is solved as a 4-th order ODE in the spectral domain $k$. For this to work, we need four initial conditions, denoted by $\hat{f}_j(k)$ for $j=0,1,2,3$. These conditions are obtained from (\ref{eq212}) and can be expressed in terms of the given amplitude $A_1(k)$ as
\begin{align}
\partial_{t}^{(j)}\hat{E}(k,0)=\hat{f}_j(k)&=(-i\omega(k))^jA_1(k)+(i\omega^*(-k))^jA_1^*(-k),\label{eq222}
\end{align}
\noindent
By choosing $A_1(k),k>0$, all four initial conditions are defined. From the relation  (\ref{eq47.1}) between $E(z,t)$ and $A(z,t)$, we then calculate the initial condition for the amplitude equation (\ref{eq40}) in exactly the same way as in (\ref{eq103}). The same procedure also applies in terms of the extra mode $B(z,t)$ whether or not it should be included. As we will see, this extra mode will be included in the following numerical tests.

We proceed to the implementation part and choose the initial amplitude $A_1(k)$ to be the same as for the toy model. A Gaussian centered at the wave number $k_0$.
\begin{align}
A_1(k)=\left\{
\begin{array}{cc}
De^{-\delta(k-k_0)^2},& k>0\\
0,& k<0 
\end{array}\right.,\label{eq224}
\end{align}
for some $D,\delta>0$. We will use different values for $k_0$ in the next two numerical tests.

\paragraph{Lorentz test for ultraviolet resonance}
The model in chapter 3.1.2 resulted in an ill-posed amplitude equation that turned out to give us a correct solution within the asymptotic regime for our perturbation scheme. In this test we picked the parameters in such a way that it results in a well-posed amplitude equation.
\begin{figure}[H]
  \centering
\captionsetup{width=0.85\textwidth}
\begin{subfigure}{.5\textwidth}
  \centering
  \includegraphics[scale=0.25]{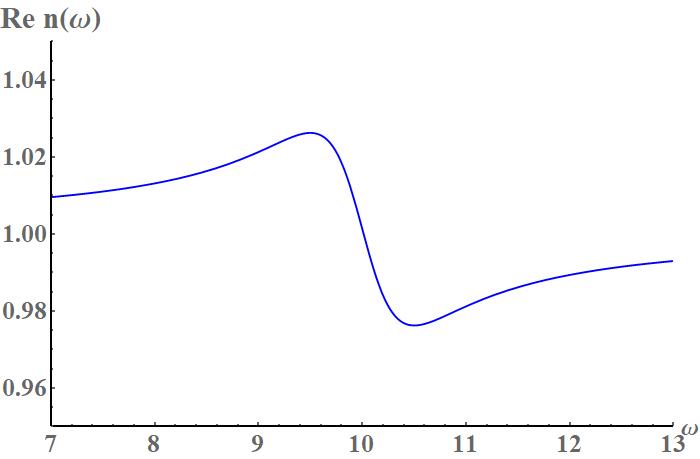}
  \caption{}
  \label{fig4a}
\end{subfigure}%
\begin{subfigure}{.5\textwidth}
  \centering
  \includegraphics[scale=0.25]{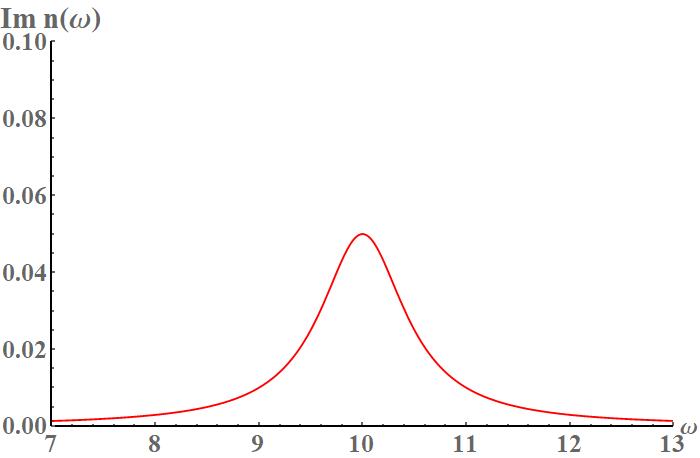}
  \caption{}
  \label{fig4b}
\end{subfigure}
\caption{Refractive index $n(\omega)$ for Lorentz test for ultraviolet resonance after scaling. The parameters in this figure are $a=-1,b=-1$ and $c=100$.}
\label{fig4}
\end{figure}
\noindent
In figure \ref{fig4} we see the real and imaginary part of the refractive index. The choice of parameters in this figure is
\begin{align}
a&=-1,\nonumber\\
b&=-1,\nonumber\\
c&=100,\label{eq224.1}
\end{align}
The parameter for the nonlinear term is again chosen to be $\varepsilon=10^{-1}$.
\begin{figure}[H]
  \centering
\captionsetup{width=0.85\textwidth}
\begin{subfigure}{.5\textwidth}
  \centering
  \includegraphics[scale=0.33]{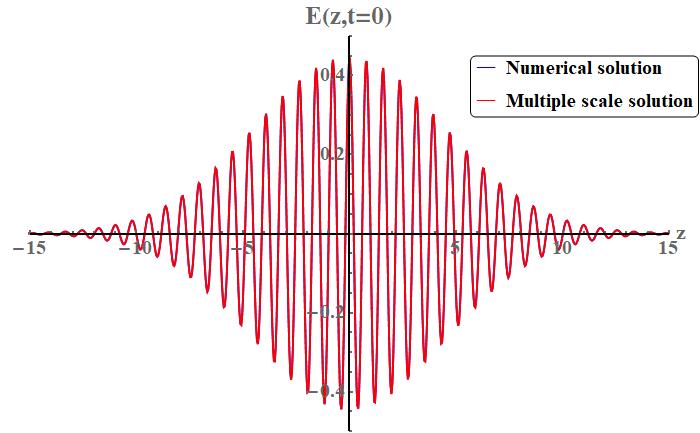}
  \caption{}
  \label{fig5a}
\end{subfigure}%
\begin{subfigure}{.5\textwidth}
  \centering
  \includegraphics[scale=0.33]{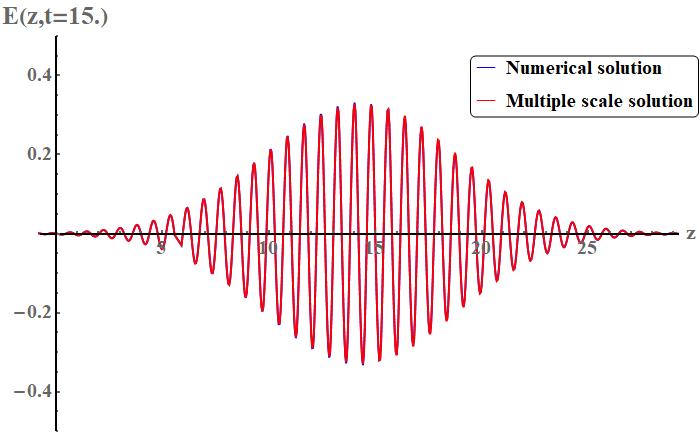}
  \caption{}
  \label{fig5b}
\end{subfigure}
\begin{subfigure}{.5\textwidth}
  \centering
  \includegraphics[scale=0.33]{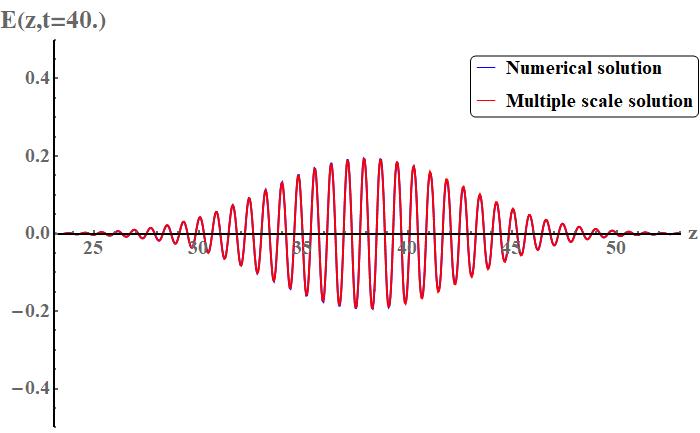}
  \caption{}
  \label{fig5c}
\end{subfigure}%
\begin{subfigure}{.5\textwidth}
  \centering
  \includegraphics[scale=0.33]{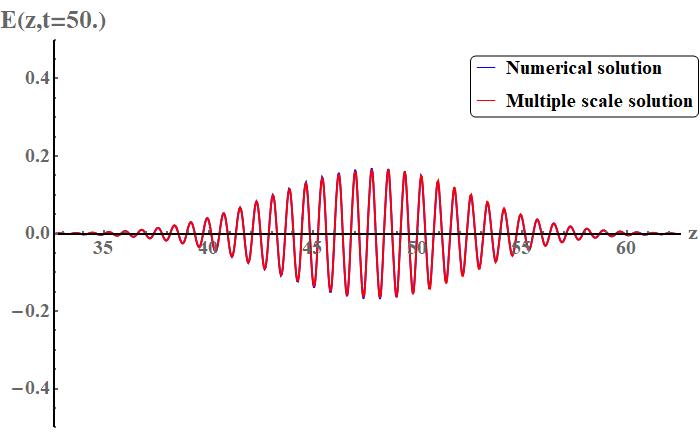}
  \caption{}
  \label{fig5d}
\end{subfigure}%
\caption{Lorentz test for ultraviolet resonance. The solution of (\ref{eq218}) using a system of first order ODE's (blue function) compared with the MMS solution (\ref{eq47.1}) (red function).}
\label{fig5}
\end{figure}
\noindent
When comparing the two solutions from both equations (\ref{eq218}) and the MMS solution  in figure \ref{fig5}, we see that they overlap sufficiently at all the presented times. However, comparing the solutions in the spectral domain can reveal the more subtle differences.

\begin{figure}[H]
  \centering
\captionsetup{width=0.85\textwidth}
\begin{subfigure}{.5\textwidth}
  \centering
  \includegraphics[scale=0.25]{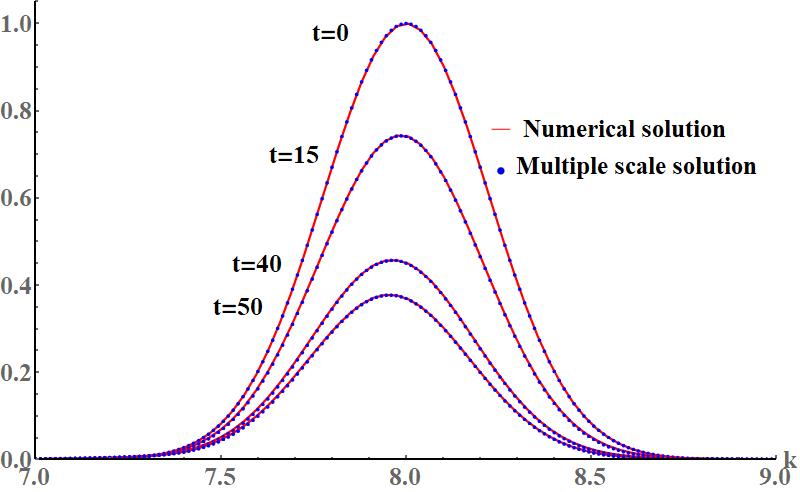}
  \caption{}
  \label{fig5.2a}
\end{subfigure}%
\begin{subfigure}{.5\textwidth}
  \centering
  \includegraphics[scale=0.25]{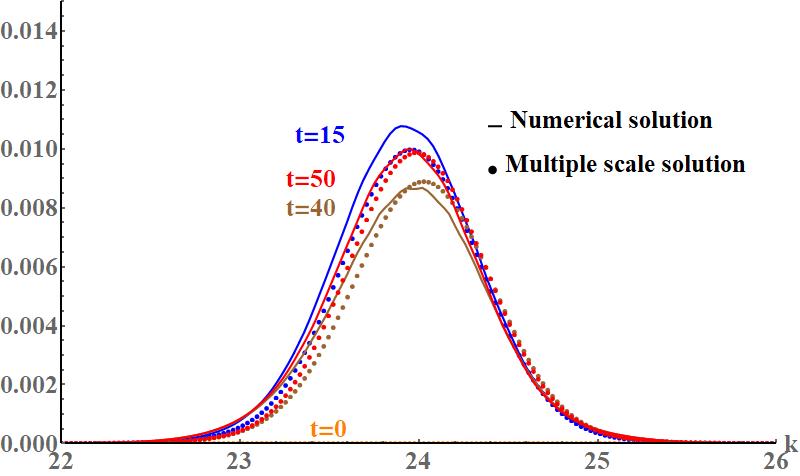}
  \caption{}
  \label{fig5.2b}
\end{subfigure}
\caption{Lorentz test for ultraviolet resonance. A close-up for the spectrum in $k$ of the solutions of (\ref{eq218}) (blue function) compared with the spectrum of the MMS solution (\ref{eq47.1}) (red function).}
\label{fig5.2}
\end{figure}
\noindent
On figure \ref{fig5.2} we see that they overlap quite nicely. The part of the spectrum we are interested in is the smaller gaussian hump centered at $3k_0=24$, where the nonlinearity manifests. In order to see the details how well these two functions overlap, we can look at the figure \ref{fig5.2} for both the gaussian humps. Figure \ref{fig5.2b} is rather convincing us of the accuracy of the MMS solution. On figure \ref{fig5.2a} is the main gaussian that represents the linear part of the solution. As we can see, the MMS solution is indeed accurate up to the order of $\varepsilon^2$ for both the toy model and the Lorentz model of dispersion. The parameters used in figures \ref{fig5} and \ref{fig5.2} are $k_0=8,\omega_0=\omega(k_0)=7.9  - i1.99\times 10^{-2}$.  In the next paragraph we run one more numerical test with a different set of parameters.

\begin{figure}[H]
  \centering
\captionsetup{width=0.85\textwidth}
  \includegraphics[scale=0.3]{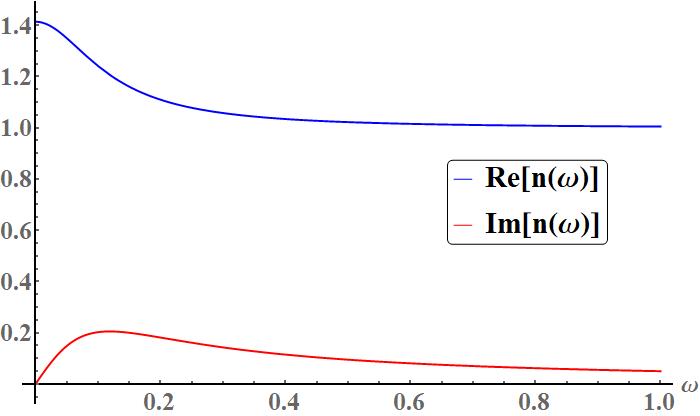}
  \caption{Refractive index $n(\omega)$ for Lorentz test for infrared resonance after scaling. The parameters in this figure are $a=-0.25,b=-10$ and $c=1$.}
\label{fig8}
\end{figure}

\paragraph{Lorentz test for infrared resonance}
Let us now choose the parameters $a,b,c$ for the Lorentz test for infrared resonance. In this model we are using the following parameters in the refracting index:
\begin{align}
a&=-0.25,\nonumber\\
b&=-10,\nonumber\\
c&=1,\label{eq226}
\end{align}
We can see the plot of this refraction index in figure (\ref{fig8}). We also notice that the graph of this function is somewhat unusual. The resonance for this function is much closer to 0 than for the index in figure \ref{fig4}. However, the resonance area in the non-scaled version of figure \ref{fig8} is located at the frequencies that are around $1.5\times 10^{14}$ (with the scaling factor $\Omega_0=1.5\times 10^{15}$) which is in the infrared range.

\noindent

The initial condition for the amplitude remains the same as in (\ref{eq224}) as well as the parameters $\delta,D,\eta$ and $E_0$. The refractive index parameters give us the scaled frequency of the pulse obtained from the dispersion relation $\omega_0=6.29 - i4.91\times 10^{-2}$.

\begin{figure}[H]
  \centering
\captionsetup{width=0.85\textwidth}
\begin{subfigure}{.5\textwidth}
  \centering
  \includegraphics[scale=0.33]{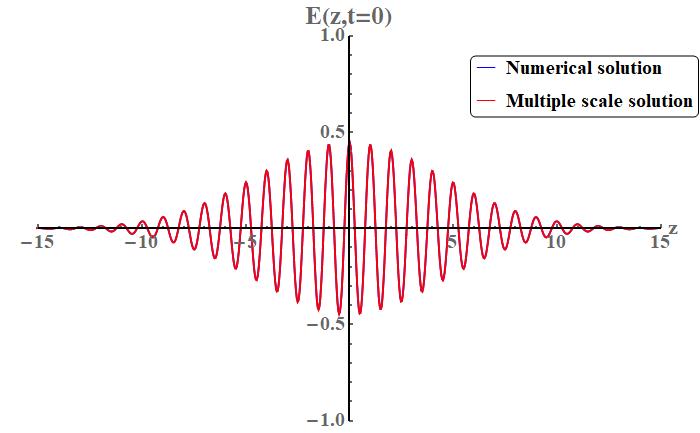}
  \caption{}
  \label{fig9a}
\end{subfigure}%
\begin{subfigure}{.5\textwidth}
  \centering
  \includegraphics[scale=0.33]{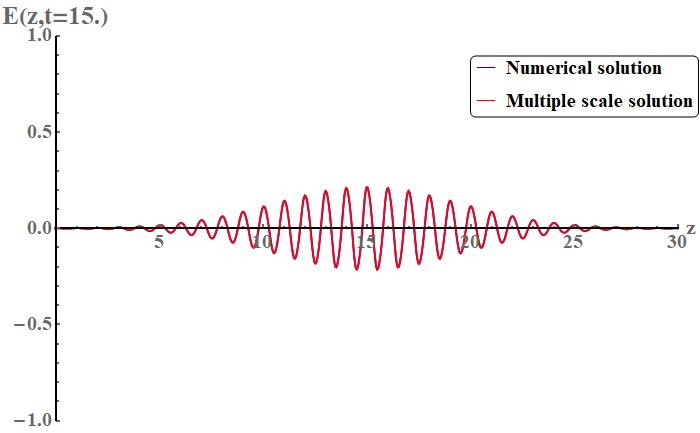}
  \caption{}
  \label{fig9b}
\end{subfigure}
\begin{subfigure}{.5\textwidth}
  \centering
  \includegraphics[scale=0.33]{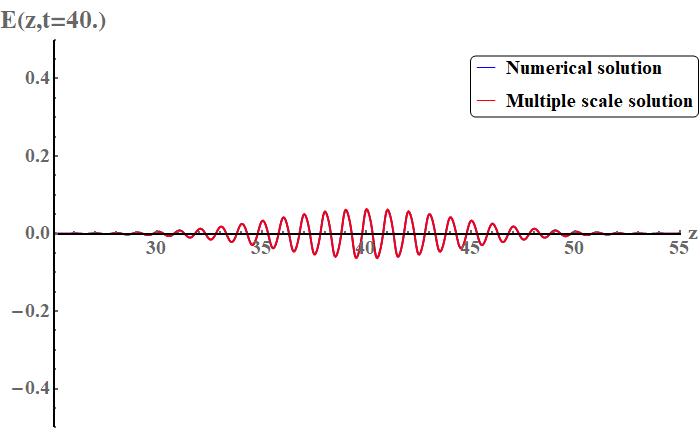}
  \caption{}
  \label{fig9c}
\end{subfigure}%
\begin{subfigure}{.5\textwidth}
  \centering
  \includegraphics[scale=0.33]{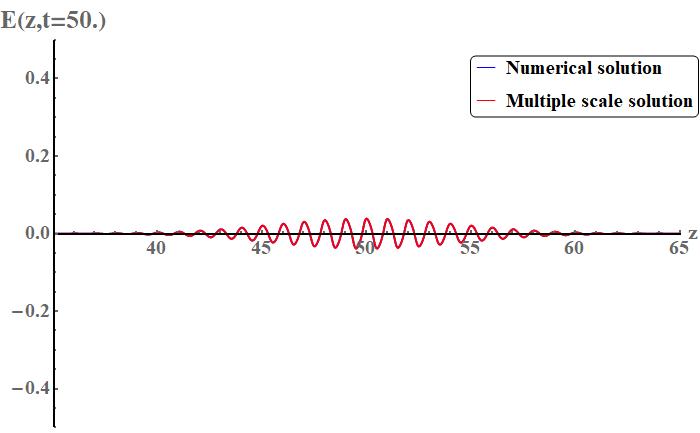}
  \caption{}
  \label{fig9d}
\end{subfigure}%
\caption{Lorentz test for infrared resonance. The solution of (\ref{eq218}) using a system of first order ODE's (blue function) compared with the MMS solution (\ref{eq47.1}) (red function).}
\label{fig9}
\end{figure}

\begin{figure}[H]
  \centering
\captionsetup{width=0.85\textwidth}
\begin{subfigure}{.5\textwidth}
  \centering
  \includegraphics[scale=0.25]{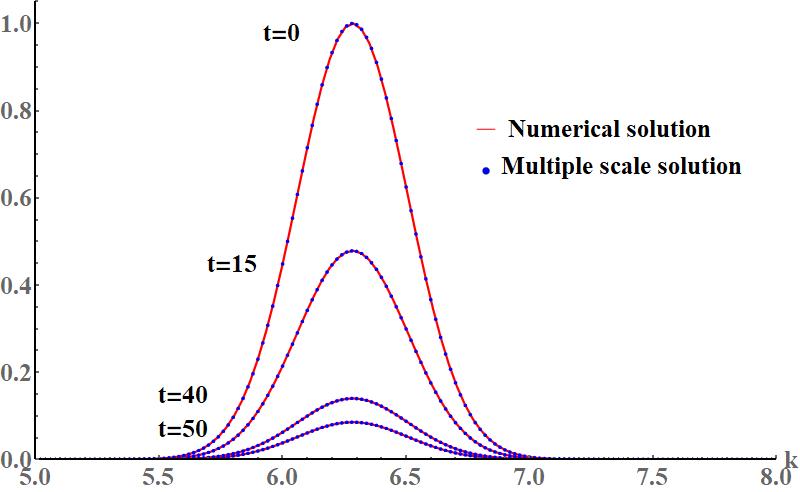}
  \caption{}
  \label{fig11a}
\end{subfigure}%
\begin{subfigure}{.5\textwidth}
  \centering
  \includegraphics[scale=0.25]{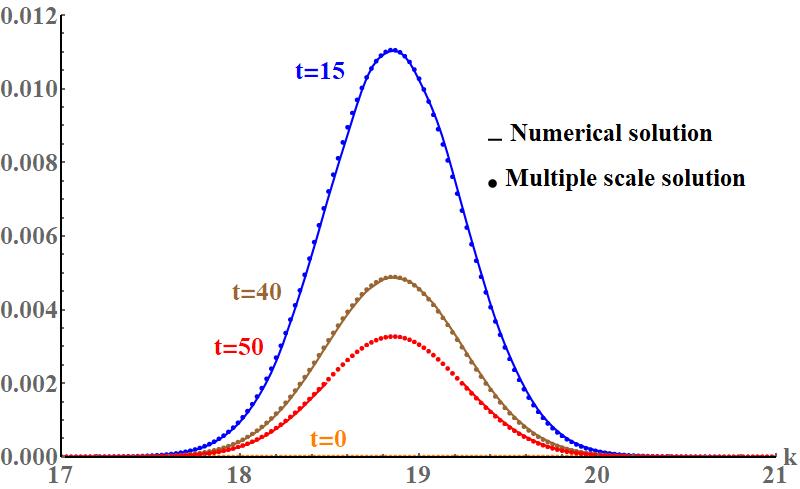}
  \caption{}
  \label{fig11b}
\end{subfigure}
\caption{Lorentz test for infrared resonance. A close-up for the spectrum in $k$ of the solutions of (\ref{eq218}) (blue function) compared with the spectrum of the MMS solution (\ref{eq47.1}) (red function).}
\label{fig11}
\end{figure}

As seen from the results in figures (\ref{fig9}) and (\ref{fig11}) we can conclude that the MMS solution proved itself and it is sufficiently accurate compared to the high precision numerical solution up to the error of order $\varepsilon^2=10^{-2}$. The parameters used in figures \ref{fig9} and \ref{fig11} are $k_0=2\pi,\omega_0=6.29 - i4.91\times 10^{-2}$.

In the next subsection we investigate the stability of the amplitude equation for both the Lorentz tests for ultraviolet and infrared resonance.

\subsubsection{Stability  and well-posedness}
\begin{figure}[t!]
  \centering
\captionsetup{width=0.85\textwidth}
\begin{subfigure}{.5\textwidth}
  \centering
  \includegraphics[scale=0.33]{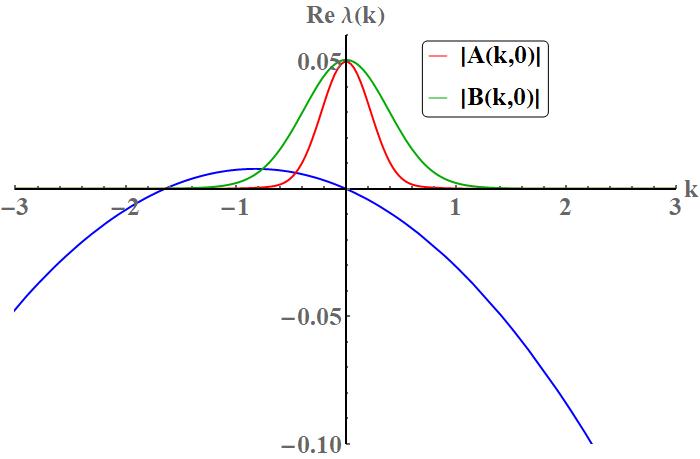}
  \caption{Lorentz test for ultraviolet resonance.}
  \label{fig6a}
\end{subfigure}%
\begin{subfigure}{.5\textwidth}
  \centering
  \includegraphics[scale=0.33]{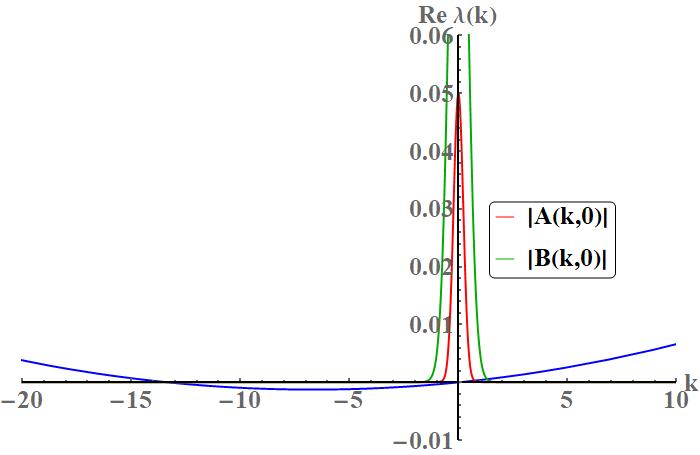}
  \caption{Lorentz test for infrared resonance.}
  \label{fig6b}
\end{subfigure}
\caption{Stability curve for the amplitude equation (\ref{eq50}). The parameters in \ref{fig6a} are $a_1=-1.15\times 10^{-2}$ and $a_2=-1.9\times 10^{-2}$ and in \ref{fig6b} are $a_1=2.83\times 10^{-5}$ and $a_2=3.77\times 10^{-4}$.}
\label{fig6}
\end{figure}

In this subsection we look at the stability of the amplitude equation (\ref{eq50}) with the parameters used in Lorentz tests for ultraviolet and infrared resonance. The formula for the stability is the same as in (\ref{eq105.2}). On figure (\ref{fig6a}) we can see the graph of the function $\text{Re}\;\lambda(k)$ for both tests. As mention earlier in Lorentz test for ultraviolet resonance, we chose the parameters for this test such that it results in a well-posed amplitude equation. Figure \ref{fig6a} confirms this. For the Lorentz test for infrared resonance we got ill-posedness again. For both cases there is a small range of $k$ where the amplitude is stable (in the ill-posed case) or unstable (in the well-posed case). The figure \ref{fig6a} tells us that we do not have to worry about the exponential growth for Lorentz model for ultraviolet resonance whereas the opposite is true for Lorentz model for infrared resonance \ref{fig6b}.

\numberwithin{equation}{section}
\section{Conclusion}
In this paper we have successfully derived the MMS solution to Maxwell's equation (\ref{eq20}) and demonstrated its numerical accuracy. During the process we introduced the Sellmeier transformation that helped us express Maxwell's equation without any pseudodifferential operators. We showed that the obtained MMS solutions provide a good approximation to the solution of the nonlinear Maxwell's equation (\ref{eq20}) up to order $\varepsilon^2$. The key features of our MMS solution are the linearity of the amplitude equation and the complex nature of the modes. The linearity made the amplitude equation analytically solvable in a much faster way than the original Maxwell's equation which is one of the main advantages of MMS. On the other hand, in some cases, the amplitude equation turned out to be ill-posed. However it does not represent a problem because of the nature of the stability curve and the location of the initial condition.

We did three numerical tests; the first corresponds to a toy model of dispersion, and the other two the, more physical, Lorentz model of dispersion. The MMS solution performed very well in all three cases. For the two Lorentz cases, we chose the frequency of the initial pulse to be in front of and behind the resonance, where the first case gave us a well-posed amplitude equation while the second one turned out to be an ill-posed one.

The key idea behind MMS is to maintain the ordering in asymptotic expansions. This ordering also applies to the Taylor expansion of the susceptibility (\ref{eq22}). It is clear, that this ordering depends on the parameters in $\hat{\chi}(\omega)$ itself, but also on the frequency $\omega_0$ around which we are expanding. In our countless numerical tests we have tried many different combinations of the parameters. We observed that this ordering failed for certain values for the parameters. Especially when one wants to be near or at resonance. Certainly, further investigation is needed to pinpoint exactly where it can go wrong. The derivatives of $\hat{\chi}(\omega)$ also appear in the constant $\alpha$, which is required to be of a certain order since it appears in the $\varepsilon^2$ part of the amplitude equation. Alternatively, different methods may be used to expand the susceptibility, for example rational function expansion.

Another problem with our approach that also relates to the values of parameters in $\hat{\chi}(\omega)$, is the growth of the nonlinear term in the Maxwell's equation. As expected, if the height of the spectral peak at $3k_0$ exceeds the order of $\varepsilon^2$ in the numerical solution, our MMS solution fails to maintain its accuracy. This problem, however, occurred less often than the aforementioned violation of ordering in the $\hat{\chi}(\omega)$ expansion. In order to resolve these issues, more extensive work is required.

\begin{appendices}
\numberwithin{equation}{section}
\section*{Appendix A}
\renewcommand{\theequation}{A.\arabic{equation}}
In this appendix we introduce a more convenient representation of the linear polarization (\ref{eq4}). First, let us mention that our convention for Fourier transform is
\begin{align}
F(\omega)&=\frac{1}{\sqrt{2\pi}}\int_{-\infty}^\infty \mathrm{d}x f(t)e^{-i\omega t},\nonumber\\
f(t)&=\frac{1}{\sqrt{2\pi}}\int_{-\infty}^\infty \mathrm{d}x F(\omega)e^{i\omega t}.\label{eq5.1}
\end{align}
Then using (\ref{eq5.1}) we have according to (\ref{eq4})
\begin{align}
\textbf{P}_L&=\varepsilon_0\int_{-\infty}^t\mathrm{d}t'\chi(t-t')\textbf{E}(\textbf{x},t'),\nonumber\\
&=\varepsilon_0\int_{-\infty}^\infty\mathrm{d}\omega\hat{\chi}(\omega)\hat{\textbf{E}}(\textbf{x},\omega)e^{-i\omega t},\nonumber\\
&=\varepsilon_0\int_{-\infty}^\infty\mathrm{d}\omega\left(\sum_{n=0}^\infty\frac{\hat{\chi}(0)}{n!}\omega^n\right)\hat{\textbf{E}}(\textbf{x},\omega)e^{-i\omega t},\nonumber\\
&=\varepsilon_0\sum_{n=0}^\infty\frac{\hat{\chi}(0)}{n!}\left(\int_{-\infty}^\infty\mathrm{d}\omega\omega^n\hat{\textbf{E}}(\textbf{x},\omega)e^{-i\omega t}\right),\nonumber\\
&=\varepsilon_0\sum_{n=0}^\infty\frac{\hat{\chi}(0)}{n!}\left(\int_{-\infty}^\infty\mathrm{d}\omega\left(i\partial_t\right)^n\hat{\textbf{E}}(\textbf{x},\omega)e^{-i\omega t}\right),\nonumber\\
&=\varepsilon_0\sum_{n=0}^\infty\frac{\hat{\chi}(0)}{n!}\left(i\partial_t\right)^n\left(\int_{-\infty}^\infty\mathrm{d}\omega\left(i\partial_t\right)^n\hat{\textbf{E}}(\textbf{x},\omega)e^{-i\omega t}\right),\nonumber\\
&=\varepsilon_0\sqrt{2\pi}\hat{\chi}\left(i\partial_t\right)\textbf{E}(\textbf{x},t),\label{eq6}
\end{align}
where $\hat{\chi}(\omega)$ is the Fourier transform of $\chi(t)$.

\setcounter{equation}{0}
\section*{Appendix B}
\renewcommand{\theequation}{B.\arabic{equation}}
In this appendix we derive the amplitude equation and the MMS solution to (\ref{eq55}). The multiple scales were introduced in (\ref{eq21}). Using these on the equation (\ref{eq55}), we get the following perturbation hierarchy
\begin{align}
\varepsilon^0:& & a\partial_{t_0t_0t_0}e_0+(\gamma+1)\partial_{t_0t_0}e_0-a\partial_{z_0z_0t_0}e_0-\gamma \partial_{z_0z_0}e_0&=0,\label{eq57}\\
\varepsilon^1:& & a\partial_{t_0t_0t_0}e_1+(\gamma+1)\partial_{t_0t_0}e_1-a\partial_{z_0z_0t_0}e_1-\gamma\partial_{z_0z_0}e_1&=\nonumber\\
& & \left(-2 \partial_{t_0t_1}-3 a \partial_{t_0t_0t_1}+a\partial_{z_0z_0t_1} +2 a \partial_{t_0z_0z_1}-2 \gamma  \partial_{t_0t_1}+2 \gamma  \partial_{z_0z_1}\right)e_0,&\label{eq58}\\
\varepsilon^2:& & a\partial_{t_0t_0t_0}e_2+(\gamma+1)\partial_{t_0t_0}e_2-a\partial_{z_0z_0t_0}e_2-\gamma\partial_{z_0z_0}e_2&=\nonumber\\
& & \left(-3 a \partial_{t_0t_0t_2}-3 a\partial_{t_0t_1t_1}+2 a\partial_{t_0z_0z_2}+a\partial_{t_0z_1z_1}+2 a\partial_{z_0z_1t_1}+a \partial_{z_0z_0t_2}\right.\nonumber\\
& & \left.-2\partial_{t_0t_2}(\gamma+1)-\partial_{t_1t_1}(\gamma+1)+2 \gamma\partial_{z_0z_2}+\gamma\partial_{z_1z_1}\right)e_0+\left(-3 a\partial_{t_0t_0t_1}+2 a\partial_{t_0z_0z_1}\right.\nonumber\\
& &\left.+a\partial_{z_0z_0t_1}-2\partial_{t_0t_1}(\gamma+1)+2 \gamma\partial_{z_0z_1}\right)e_1-\left(a\partial_{t_0t_0t_0} +\gamma\partial_{t_0t_0}\right)e_0^3&.\label{eq59}
\end{align}
For the $\varepsilon^0$ order equation we choose the wave packet solution
\begin{align}
e_0(z_0,t_0,z_1,t_1,\ldots)=A_0(z_1,t_1,\ldots)e^{i\theta_0}+(*),\label{eq60}
\end{align}
where
\begin{align}
\theta_0=kz_0-\omega t_0,\label{eq61}
\end{align}
and where $\omega=\omega(k)$ is a complex number and a solution to the dispersion relation
\begin{align}
a\omega^3+i\omega^2(\gamma+1)-ak^2\omega-i\gamma k^2=0.\label{eq62}
\end{align}
Note that $k$ is the scaled initial spatial frequency and is used to obtain the complex frequency $\omega$.

We now proceed to the $\varepsilon$ order equation. Inserting (\ref{eq60}) into (\ref{eq58}) we get
\begin{align}
&a\partial_{t_0t_0t_0}e_1+(\gamma+1)\partial_{t_0t_0}e_1-a\partial_{z_0z_0t_0}e_1-\gamma\partial_{z_0z_0}e_1=\nonumber\\
&\left(2i\omega \partial_{t_1}(\gamma+1)A_0+3 a\omega^2 \partial_{t_1}A_0-ak^2 \partial_{t_1}A_0 +2 a\omega k \partial_{z_1}A_0+2ik \gamma  \partial_{z_1}A_0\right)e^{i\theta_0}+(*).\label{eq63}
\end{align}
In order to remove secular terms we postulate that
\begin{align}
2i\omega \partial_{t_1}(\gamma+1)A_0+3 a\omega^2 \partial_{t_1}A_0-ak^2 \partial_{t_1}A_0 +2 a\omega k \partial_{z_1}A_0+2ik \gamma  \partial_{z_1}A_0&=0,\nonumber\\
\partial_{t_1}A_0\left(2i\omega(\gamma+1)+3a\omega^2-ak^2\right)+\partial_{z_1}A_0\left(2ak\omega+2i\gamma k\right)&=0.\label{eq64}
\end{align}
Observe that by differentiating the dispersion relation (\ref{eq62}) with respect to $k$ we get
\begin{align}
a\omega^3+i\omega^2(\gamma+1)-ak^2\omega-i\gamma k^2&=0,\nonumber\\
&\Downarrow\nonumber\\
3a\omega^2\omega'(k)+i2\omega\omega'(k)(\gamma+1)-a2k\omega-ak^2\omega'(k)-i2\gamma k&=0,\nonumber\\
3a\omega^2+i2\omega(\gamma+1)-ak^2&=\frac{a2k\omega+i2\gamma k}{\omega'(k)},\label{eq65}
\end{align}
Using (\ref{eq65}), the equation (\ref{eq64}) can be written in the form
\begin{align}
\partial_{t_1}A_0+\omega'(k)\partial_{z_1}A_0=0.\label{eq66}
\end{align}
The $\varepsilon$ order equation (\ref{eq63}) now simplifies into
\begin{align}
a\partial_{t_0t_0t_0}e_1+(\gamma+1)\partial_{t_0t_0}e_1-a\partial_{z_0z_0t_0}e_1-\gamma\partial_{z_0z_0}e_1=0.\label{eq67}
\end{align}
We choose the special zero solution to (\ref{eq67})
\begin{align}
e_1=0.\label{eq68}
\end{align}
We now compute the right-hand side of the order $\varepsilon^2$ equation. Inserting (\ref{eq68}) into the right-hand side of the order $\varepsilon^2$ equation (\ref{eq59}) we get
\begin{align}
a\partial_{t_0t_0t_0}e_2+(\gamma+1)\partial_{t_0t_0}e_2-a\partial_{z_0z_0t_0}e_2-\gamma\partial_{z_0z_0}e_2&=\left( 3a\omega^2\partial_{t_2}A_0+3ia\omega\partial_{t_1t_1}A_0+2ak\omega\partial_{z_2}A_0\right.\nonumber\\
&-ia\omega\partial_{z_1z_1}A_0+2iak\partial_{t_1z_1}A_0-ak^2\partial_{t_2}A_0\nonumber\\
&+2i\omega(\gamma+1)\partial_{t_2}A_0-(\gamma+1)\partial_{t_1t_1}A_0\nonumber\\
&\left.+2i\gamma k\partial_{z_2}A_0+\gamma n_0^2\partial_{z_1z_1}A_0\right)e^{i\theta_0}-NST+(*),\label{eq69}
\end{align}
where
\begin{align}
NST&=a27i\omega^3A_0^3 e^{3i\theta_0}+A_0^2A_0^*e^{i\theta_0}e^{2t_0\omega_i}\left( 3ia\omega^3-18a\omega^2\omega_i-36ia\omega\omega_i^2 +24a\omega_i^3 \right)\nonumber\\
& -9\gamma\omega^2A_0^2e^{3i\theta_0}+3\gamma A_0^2A_0^*e^{i\theta_0}e^{2t_0\omega_i}\left(2\omega_i-i\omega\right)^2,\label{eq70}
\end{align}
are the non-secular terms and where $\omega_i=\text{Im}\;\omega$. In order to remove secular terms we postulate that
\begin{align}
&3a\omega^2\partial_{t_2}A_0+3ia\omega\partial_{t_1t_1}A_0+2ak\omega\partial_{z_2}A_0-ia\omega\partial_{z_1z_1}A_0+2iak\partial_{t_1z_1}A_0-ak^2\partial_{t_2}A_0\nonumber\\
&+2i\omega(\gamma+1)\partial_{t_2}A_0-(\gamma+1)\partial_{t_1t_1}A_0+2i\gamma k\partial_{z_2}A_0+\gamma\partial_{z_1z_1}A_0=0,\nonumber\\
&\quad\quad\quad\quad\quad\quad\quad\quad\quad\quad\quad\quad\quad\quad\quad\quad\Downarrow\nonumber\\
&\partial_{t_2}A_0\left(3a\omega^2-ak^2+2i\omega(\gamma+1)\right)+\partial_{t_1t_1}A_0\left(3ia\omega-(\gamma+1)\right)+\partial_{z_2}A_0\left(2ak\omega+2i\gamma k\right)\nonumber\\
&+2iak\partial_{t_1z_1}A_0+\partial_{z_1z_1}A_0\left(\gamma-ia\omega\right)=0.\label{eq71}
\end{align}
From the equation (\ref{eq64}) we express the $\partial_{t_1z_1}$ derivative in terms of $\partial_{t_1t_1}$ in the following way
\begin{align}
\partial_{t_1}A_0\left(2i\omega(\gamma+1)+3a\omega^2-ak^2\right)&=-\partial_{z_1}A_0\left(2ak\omega+2i\gamma k\right)\quad\quad/\partial_{t_1},\nonumber\\
\partial_{t_1t_1}A_0\frac{ak^2-2i\omega(\gamma+1)-3a\omega^2}{2ak\omega+2i\gamma k}&=\partial_{t_1z_1}A_0.\label{eq72}
\end{align}
Substituting (\ref{eq72}) back to (\ref{eq71}) and using the relation (\ref{eq65}) by then term $\partial_{t_2}A_0$ we get
\begin{align}
&\partial_{t_2}A_0\frac{a2k\omega+i2\gamma k}{\omega'(k)}+\partial_{t_1t_1}A_0\left(3ia\omega-(\gamma+1)\right)+\partial_{z_2}A_0\left(2ak\omega+2i\gamma k\right)\nonumber\\
&+2iak\frac{ak^2-2i\omega(\gamma+1)-3a\omega^2}{2ak\omega+2i\gamma k}\partial_{t_1t_1}A_0+\partial_{z_1z_1}A_0\left(\gamma-ia\omega\right)=0,\nonumber\\
&\quad\quad\quad\quad\quad\quad\quad\quad\quad\quad\quad\quad\quad\quad\quad\quad\Downarrow\nonumber\\
&\partial_{t_2}A_0+\omega'(k)\partial_{z_2}A_0+\frac{\omega'(k)}{a2k\omega+i2\gamma k}\left(3ia\omega-(\gamma+1)+ 2iak\frac{ak^2-2i\omega(\gamma+1)-3a\omega^2}{2ak\omega+2i\gamma k}\right)\partial_{t_1t_1}A_0\nonumber\\
&-i\frac{\omega'(k)}{2k}\partial_{z_1z_1}A_0=0,\label{eq73}
\end{align} 
Let us deal with the factor by the term $\partial_{t_1t_1}A_0$ separately.
\begin{align}
&\frac{\omega'(k)}{a2k\omega+i2\gamma k}\left(3ia\omega-(\gamma+1)+ 2iak\frac{ak^2-2i\omega(\gamma+1)-3a\omega^2}{2ak\omega+2i\gamma k}\right)\nonumber\\
&=\frac{i\omega'(k)}{2k}\left(\frac{\gamma+1-3ia\omega}{\gamma-ia\omega}+\frac{a\left(2i\omega(\gamma+1)+3a\omega^2-ak^2\right)}{(\gamma-ia\omega)^2}\right)\nonumber\\
&=\frac{i\omega'(k)}{2k}\frac{\left(\gamma+1-3ia\omega\right)(\gamma-ia\omega)+a2i\omega(\gamma+1)+3a^2\omega^2-a^2k^2}{(\gamma-ia\omega)^2}\nonumber\\
&=\frac{i\omega'(k)}{2k}\frac{\gamma^2+\gamma-2ia\omega\gamma+ia\omega-a^2k^2}{(\gamma-ia\omega)^2}.\label{eq74}
\end{align}
Now we turn once again for help to the dispersion relation (\ref{eq62}) and find
\begin{align}
a\omega^3+i\omega^2(\gamma+1)-ak^2\omega-i\gamma k^2=0,\nonumber\\
&\Downarrow\nonumber\\
i\omega^2(\gamma-ia\omega)+i\omega^2-ik^2(\gamma-ia\omega)&=0,\nonumber\\
k^2&=\frac{\omega^2}{\gamma-ia\omega}+\omega^2.\label{eq75}
\end{align}
Inserting (\ref{eq75}) with $k=k$ back into (\ref{eq74}) we obtain
\begin{align}
\frac{i\omega'(k)}{2k}\frac{\gamma^2+\gamma-2ia\omega\gamma+ia\omega-a^2k^2}{(\gamma-ia\omega)^2}&=\frac{i\omega'(k)}{2k}\frac{\gamma^2+\gamma-2ia\omega\gamma+ia\omega-a^2\frac{\omega^2}{\gamma-ia\omega}-a^2\omega^2}{(\gamma-ia\omega)^2}\nonumber\\
&=\frac{i\omega'(k)}{2k}\frac{(\gamma-ia\omega)^2+\gamma+ia\omega-a^2\frac{\omega^2}{\gamma-ia\omega}}{(\gamma-ia\omega)^2}\nonumber\\
&=\frac{i\omega'(k)}{2k}\frac{(\gamma-ia\omega)^2+\frac{(\gamma+ia\omega)(\gamma-ia\omega)-a^2\omega^2}{\gamma-ia\omega}}{(\gamma-ia\omega)^2}\nonumber\\
&=\frac{i\omega'(k)}{2k}\left(1+\frac{\frac{\gamma^2+a^2\omega^2-a^2\omega^2}{\gamma-ia\omega}}{(\gamma-ia\omega)^2}\right)=\frac{i\omega'(k)}{2k}\left(1+\frac{\gamma^2}{(\gamma-ia\omega)^3}\right).\label{eq76}
\end{align}
Using this result in (\ref{eq73}), our final amplitude equation is
\begin{align}
\partial_{t_2}A_0+\omega'(k)\partial_{z_2}A_0+\frac{i\omega'(k)}{2k}\left(1+\frac{\gamma^2}{(\gamma-ia\omega)^3}\right)\partial_{t_1t_1}A_0-i\frac{\omega'(k)}{2k}\partial_{z_1z_1}A_0=0.\label{eq77}
\end{align}
By removing the secular terms from the equation (\ref{eq69}), the order $\varepsilon^2$ equation turns to
\begin{align}
&a\partial_{t_0t_0t_0}e_2+(\gamma+1)\partial_{t_0t_0}e_2-a\partial_{z_0z_0t_0}e_2-\gamma\partial_{z_0z_0}e_2=-\left(a27i\omega^3A_0^3 e^{3i\theta_0}\right.\nonumber\\
&\left.+A_0^2A_0^*e^{i\theta_0}e^{2t_0\omega_i}\left( 3ia\omega^3-18a\omega^2\omega_i-36ia\omega\omega_i^2 +24a\omega_i^3 \right)-9\gamma\omega^2A_0^32e^{3i\theta_0}\right.\nonumber\\
&\left. +3\gamma A_0^2A_0^*e^{i\theta_0}e^{2t_0\omega_i}\left(2\omega_i-i\omega\right)^2\right)+(*),\label{eq78}
\end{align}
which we solve for $e_2$ taking only the particular solution.
\begin{align}
e_2(z_0,t_0,\ldots)=c_1A_0^3e^{i3\theta_0}+c_2|A_0|^2A_0e^{i\theta_0}e^{2\omega_i}+(*),\label{eq79}
\end{align}
where
\begin{align}
c_1&=\frac{9\omega^2(\gamma-3ia\omega)}{-9\omega^2(\gamma-3ia\omega)+9k^2(\gamma-3ia\omega)-9\omega^2}=\frac{1}{-1+n^2(\omega)-1/(\gamma-i3a\omega)}\nonumber\\
&=\frac{1}{n^2(\omega)-n^2(3\omega)},\label{eq80}\\
c_2&=\frac{-\left(3a(2\omega_{0i}-i\omega)^3+3\gamma(2\omega_{0i}-i\omega)^2\right)}{a(2\omega_{0i}-2i\omega)^3+k^2(\gamma+a(2\omega_{0i}-i\omega))+(2\omega_{0i}-i\omega)^2(\gamma+1)}\nonumber\\
&=\frac{-3(2\omega_{0i}-i\omega)^2}{k^2+(2\omega_{0i}-i\omega)^2\left(1+1/(\gamma-ia(\omega+2i\omega_{0i}))\right)}\nonumber\\
&=\frac{3(\omega+2i\omega_{0i})^2}{k^2-(\omega+2i\omega_{0i})^2\left(1+1/(\gamma-ia(\omega+2i\omega_{0i}))\right)}\nonumber\\
&=\frac{3(\omega+2i\omega_i)^2}{k^2-\left(1+\hat{\chi}(\omega+i2\omega_i)\right) \left(\omega +2 i \omega_i\right) ^2},\label{eq81}
\end{align}
where $\omega_{i}=\text{Im}\;\omega$.

Defining as before the amplitude as in (\ref{eq45}) and proceeding the usual way using (\ref{eq66}) and (\ref{eq77}) we get the amplitude equation
\begin{align}
\partial_tA+\omega'(k)\partial_zA-i\beta\partial_{zz}A+i\alpha\partial_{tt}A=0,\label{eq82}
\end{align}
where
\begin{align}
\alpha&=\frac{\omega'(k)}{2k}\left(1+\frac{\gamma^2}{(\gamma-ia\omega)^3}\right),\label{eq83}\\
\beta&=\frac{\omega'(k)}{2k}.\label{eq84}
\end{align}
The overall approximate solution to (\ref{eq55}) is then
\begin{align}
E(z,t)=A(z,t)e^{i(kz-\omega t)}+c_1\varepsilon^2A^3(z,t)e^{i3(kz-\omega t)}+c_2\varepsilon^2|A(z,t)|^2A(z,t)e^{i(kz-\omega t)}e^{2t\omega_i}+(*),\label{eq85}
\end{align}
where $c_1,c_2$ are defined in (\ref{eq80}) and (\ref{eq81}). To verify that $\alpha,\beta$ in (\ref{eq47}) and (\ref{eq82}) are the same, we look at (\ref{eq39}) and get
\begin{align}
&\omega'(k)\frac{n^2(\omega)+2\omega a\sqrt{2\pi}\hat{\chi}'(\omega)+\frac{a^2}{2}\omega^2\sqrt{2\pi}\hat{\chi}''(\omega)}{2k}=\frac{\omega'(k)}{2k}\left(1+\frac{1}{\gamma-ia\omega}+2\omega a\frac{i}{(\gamma-ia\omega)^2}\right.\nonumber\\
&\left.+\frac{a^2}{2}\omega^2\frac{-2}{(\gamma-ia\omega)^3}\right)=\frac{\omega'(k)}{2k}\left(\frac{(\gamma-ia\omega)^3+(\gamma-ia\omega)^2+2\omega ai(\gamma-ia\omega)-a^2\omega^2}{(\gamma-ia\omega)^3}\right)\nonumber\\
&=\frac{\omega'(k)}{2k}\left(\frac{(\gamma-ia\omega)^3+\gamma^2}{(\gamma-ia\omega)^3}\right)=\frac{\omega'(k)}{2k}\left(1+\frac{\gamma^2}{(\gamma-ia\omega)^3}\right),\label{eq86}
\end{align}
which is the same as (\ref{eq83}).

Using the same argument as in (\ref{eq48}) we can simplify the amplitude equation (\ref{eq82}) into
\begin{align}
\partial_tA+\omega'(k)\partial_zA-i\partial_{zz}A\left(\beta-\alpha\left(\omega'(k)\right)^2\right)=0.\label{eq87}
\end{align}

\setcounter{equation}{0}
\section*{Appendix C}
\renewcommand{\theequation}{C.\arabic{equation}}
We are going to use multiple scale method to derive an approximate solution to (\ref{eq164}). Using (\ref{eq21}) on the equation (\ref{eq164}) we get the perturbation hierarchy
\begin{align}
\varepsilon^0:& &(1+c)\partial_{t_0t_0}e_0-b\partial_{t_0t_0t_0}e_0-a\partial_{t_0t_0t_0t_0}e_0-c\partial_{z_0z_0}e_0+b\partial_{t_0z_0z_0}e_0+a\partial_{t_0t_0z_0z_0}e_0=0,&\label{eq165}\\
\varepsilon^1:& &(1+c)\partial_{t_0t_0}e_1-b\partial_{t_0t_0t_0}e_1-a\partial_{t_0t_0t_0t_0}e_1-c\partial_{z_0z_0}e_1+b\partial_{t_0z_0z_0}e_1+a\partial_{t_0t_0z_0z_0}e_1=&\nonumber\\
& & -2(1+c)\partial_{t_0t_1}e_0+3 b \partial_{t_0t_0}\partial_{t_1}e_0 + 4 a \partial_{t_0t_0t_0}\partial_{t_1}e_0 - b \partial_{t_1} \partial_{z_0z_0}e_0 - 
 2 a \partial_{t_0} \partial_{t_1} \partial_{z_0z_0}e_0 + 2 c \partial_{z_0} \partial_{z_1}e_0&\nonumber\\
 & & - 2 b \partial_{t_0} \partial_{z_0} \partial_{z_1}e_0 - 
 2 a \partial_{t_0t_0}\partial_{z_0}\partial_{z_1}e_0,&\label{eq166}\\
\varepsilon^2:& &(1+c)\partial_{t_0t_0}e_2-b\partial_{t_0t_0t_0}e_2-a\partial_{t_0t_0t_0t_0}e_2-c\partial_{z_0z_0}e_2+b\partial_{t_0z_0z_0}e_2+a\partial_{t_0t_0z_0z_0}e_2=&\nonumber\\
& &-(1 +c)\partial_{t_1t_1}e_0 + 3 b \partial_{t_0}\partial_{t_1t_1}e_0 + 
 6 a\partial_{t_0t_0}\partial_{t_1t_1}e_0 + -2(1+c)\partial_{t_0}\partial_{t_2}e_0 + 3 b \partial_{t_0t_0}\partial_{t_2}e_0&\nonumber\\
 & & +4 a \partial_{t_0t_0t_0}\partial_{t_2}e_0 - a \partial_{t_1t_1}\partial_{z_0z_0}e_0 - b \partial_{t_2} \partial_{z_0z_0}e_0 - 
 2 a \partial_{t_0} \partial_{t_2} \partial_{z_0z_0}e_0- 2 b \partial_{t_1} \partial_{z_0} \partial_{z_1}e_0&\nonumber\\
 & &  - 4 a \partial_{t_0} \partial_{t_1} \partial_{z_0} \partial_{z_1}e_0 + 
 c \partial_{z_1z_1}e_0 - b \partial_{t_0} \partial_{z_1z_1}e_0 - a \partial_{t_0t_0} \partial_{z_1z_1}e_0 + 2 c \partial_{z_0} \partial_{z_2}e_0- 
 2 b \partial_{t_0} \partial_{z_0} \partial_{z_2}e_0&\nonumber\\
 & &  - 2 a \partial_{t_0t_0}\partial_{z_0}\partial_{z_2}e_0-2(1 +c0) \partial_{t_0}\partial_{t_1}e_1 + 3 b \partial_{t_0t_0}\partial_{t_1}e_1 + 4 a \partial_{t_0t_0t_0} \partial_{t_1}e_1 -  b \partial_{t_1}\partial_{z_0z_0}e_1&\nonumber\\
 & & - 2 a \partial_{t_0} \partial_{t_1}\partial_{z_0z_0}e_1 + 2 c \partial_{z_0}\partial_{z_1}e_1 - 
 2 b \partial_{t_0}\partial_{z_0}\partial_{z_1}e_1 - 2 a\partial_{t_0t_0}\partial_{z_0}\partial_{z_1}e_1&\nonumber\\
 & & -\varepsilon^2 (c \partial_{t_0t_0} - b \partial_{t_0t_0t_0} - 
   a \partial_{t_0t_0t_0t_0})e_0^3.&\label{eq167}
\end{align}
The order $\varepsilon^0$ equation has a solution
\begin{align}
e_0(z_0,t_0,\ldots)=A_0(z_1,t_1,\ldots)e^{i\theta_0}+(*),\label{eq168}
\end{align}
where
\begin{align}
\theta_0=kz_0-\omega t_0,\label{eq169}
\end{align}
and where $\omega=\omega(k)$ a complex number and a solution to the dispersion relation
\begin{align}
&a \omega ^4+i b \omega ^3+\omega ^2(c+1)=k^2 \left(a \omega ^2+i b \omega +c\right).\label{eq170}
\end{align}
We are going to use (\ref{eq168}) in the right hand side of the equation (\ref{eq166}) and obtain
\begin{align}
&(1+c)\partial_{t_0t_0}e_1-b\partial_{t_0t_0t_0}e_1-a\partial_{t_0t_0t_0t_0}e_1-c\partial_{z_0z_0}e_1+b\partial_{t_0z_0z_0}e_1+a\partial_{t_0t_0z_0z_0}e_1=\nonumber\\
&\left(\partial_{t_1}A_0 \left(-2 i a k^2 \omega +4 i a \omega ^3+b k^2-3 b \omega ^2-i (-2 c-2) \omega \right)\right.\nonumber\\
&\left.+\partial_{z_1}A_0 \left(2 i a k \omega ^2-2 b k \omega +2 i c k\right)\right)e^{i\theta_0}+(*).\label{eq171}
\end{align}
In order to remove secular terms we postulate that
\begin{align}
&\partial_{t_1}A_0 \left(-2 i a k^2 \omega +4 i a \omega ^3+b k^2-3 b \omega ^2-i (-2 c-2) \omega \right)+\partial_{z_1}A_0 \left(2 i a k \omega ^2-2 b k \omega +2 i c k\right)=0.\label{eq172}
\end{align}
Observe that from the dispersion relation (\ref{eq62}) we have by differentiating it with respect to $k$
\begin{align}
a \omega ^4+i b \omega ^3+\omega ^2(c+1)&=k^2 \left(a \omega ^2+i b \omega +c\right),\nonumber\\
&\Downarrow\nonumber\\
\omega'(k)\left(4a' \omega ^3+i3 b \omega ^2+2\omega(c+1)\right)&=2k \left(a \omega ^2+i b \omega +c\right)+k^2\omega'(k) \left(2a' \omega+i b\right),\nonumber\\
\omega'(k)\left(4ia' \omega ^3-3 b \omega ^2+2i\omega(c+1)-2ik^2a' \omega+k^2 b\right)&=2ik a \omega ^2-2k b \omega +2ikc\label{eq173}
\end{align}
Thus (\ref{eq172}) can be written in the form
\begin{align}
\partial_{t_1}A_0+\omega'(k)\partial_{z_1}A_0=0,\label{eq174}
\end{align}
The order $\varepsilon^1$ equation simplifies into
\begin{align}
&(1+c)\partial_{t_0t_0}e_1-b\partial_{t_0t_0t_0}e_1-a\partial_{t_0t_0t_0t_0}e_1-c\partial_{z_0z_0}e_1+b\partial_{t_0z_0z_0}e_1+a\partial_{t_0t_0z_0z_0}e_1=0\label{eq175}
\end{align}
We choose the special zero solution to (\ref{eq175})
\begin{align}
e_1=0.\label{eq176}
\end{align}
We now compute the right-hand side of the order $\varepsilon^2$ equation. Inserting (\ref{eq168}) and (\ref{eq176}) into the right-hand side of the order $\varepsilon^2$ equation (\ref{eq167}) we get
\begin{align}
&a^2\partial_{t_0}^{(6)}e_2+2a'b\partial_{t_0}^{(5)}e_2+(b^2-2a'-2a'c)\partial_{t_0}^{(4)}e_2-(2b'+2b'c)\partial_{t_0}^{(3)}e_2+(1+c)^2\partial_{t_0}^{(2)}e_2\nonumber\\
&-a^2\partial_{t_0}^{(4)}\partial_{z_0z_0}e_2 -2a'b\partial_{t_0}^{(3)}\partial_{z_0z_0}e_2+(2a'c-b^2)\partial_{t_0}^{(2)}\partial_{z_0z_0}e_2+2b'c\partial_{t_0}^{(1)}\partial_{z_0z_0}e_2-c^2\partial_{z_0z_0}e_2=\nonumber\\
&\left(\partial_{t_1t_1}A_0 \left(a k^2-6 a \omega ^2-3 i b \omega -c-1\right)+\partial_{t_2}A_0 \left(-2 i a k^2 \omega +4 i a \omega ^3+b k^2-3 b \omega ^2-i (-2 c-2) \omega \right)\right.\nonumber\\
&\left.+\partial_{z_1z_1}A_0 \left(a \omega ^2+i b \omega +c\right)+\partial_{z_2}A_0 \left(2 i a k \omega ^2-2 b k \omega +2 i c k\right)+\partial_{t_1z_1}A_0 (-4 a k \omega -2 i b k)\right)e^{i\theta}\nonumber\\
&-NST+(*),\label{eq177}
\end{align}
where
\begin{align}
NST&=A_0^3e^{3i\theta_0}\left(-81 a \omega ^4-27 i b \omega ^3-9 c \omega ^2\right)\nonumber\\
&+3A_0^2A_0^*e^{i\theta_0}e^{2t_0\omega_i}\left(-a \omega ^4-8 i a \omega ^3 \omega_i+24 a \omega ^2 \omega_i^2+32 i a \omega  \omega_i^3-16 a \omega_i^4\right.\nonumber\\
&\left.-i b \omega ^3+6 b \omega ^2 \omega_i+12 i b \omega  \omega_i^2-8 b \omega_i^3-c \omega ^2-4 i c \omega  \omega_i+4 c \omega_i^2\right),\label{eq178}
\end{align}
are the non-secular terms and where $\omega_i=\text{Im}\;\omega$. In order to remove secular terms we postulate that
\begin{align}
&\partial_{t_1t_1}A_0 \left(a k^2-6 a \omega ^2-3 i b \omega -c-1\right)+\partial_{t_2}A_0 \left(-2 i a k^2 \omega +4 i a \omega ^3+b k^2-3 b \omega ^2-i (-2 c-2) \omega \right)\nonumber\\
&+\partial_{z_1z_1}A_0 \left(a \omega ^2+i b \omega +c\right)+\partial_{z_2}A_0 \left(2 i a k \omega ^2-2 b k \omega +2 i c k\right)+\partial_{t_1z_1}A_0 (-4 a k \omega -2 i b k)=0.\label{eq179}
\end{align}
From the equation (\ref{eq172}) we express the $\partial_{t_1z_1}$ derivative in terms of $\partial_{t_1t_1}$ in the following way
\begin{align}
&\partial_{t_1}A_0 \left(-2 i a k^2 \omega +4 i a \omega ^3+b k^2-3 b \omega ^2-i (-2 c-2) \omega \right)=-\partial_{z_1}A_0 \left(2 i a k \omega ^2-2 b k \omega +2 i c k\right),\nonumber\\
&\quad\quad\quad\quad\quad\quad\quad\quad\quad\quad\quad\quad\quad\quad\quad\quad\Downarrow\nonumber\\
&\partial_{t_1t_1}A_0 \left(2 i a k^2 \omega -4 i a \omega ^3-b k^2+3 b \omega ^2+i (-2 c-2) \omega \right)/\left(2 i a k \omega ^2-2 b k \omega +2 i c k\right)=\partial_{t_1z_1}A_0\label{eq180}
\end{align}
Substituting (\ref{eq180}) back to (\ref{eq179}) and using the relation (\ref{eq173}) on the term by $\partial_{t_2}A_0$ we get
\begin{align}
&\partial_{t_1t_1}A_0 \left(a k^2-6 a \omega ^2-3 i b \omega -c-1\right)+\partial_{t_2}A_0 \frac{1}{\omega'(k)}\left(2ik a \omega ^2-2k b \omega +2ikc \right)\nonumber\\
&+\partial_{z_1z_1}A_0 \left(a \omega ^2+i b \omega +c\right)+\partial_{z_2}A_0 \left(2 i a k \omega ^2-2 b k \omega +2 i c k\right)\nonumber\\
&+\partial_{t_1t_1}A_0 \frac{\left(2 i a k^2 \omega -4 i a \omega ^3-b k^2+3 b \omega ^2+i (-2 c-2) \omega \right)}{\left(2 i a k \omega ^2-2 b k \omega +2 i c k\right)} (-4 a k \omega -2 i b k)=0.\label{eq181}
\end{align}
We divide the whole equation by the factor in front of $\partial_{t_2}A_0$ and get
\begin{align}
&\partial_{t_2}A_0-\partial_{z_1z_1}A_0 i\frac{\omega'(k)}{2k}+\partial_{z_2}A_0 \omega'(k)+\frac{\omega'(k)}{2ik a \omega ^2-2k b \omega +2ikc}\partial_{t_1t_1}A_0 \left(\frac{p_2}{p_3}+p_1\right)=0,\label{eq182}
\end{align}
where
\begin{align}
p_1&=a k^2-6 a \omega ^2-3 i b \omega -c-1\nonumber\\
p_2&=\left(2 i a k^2 \omega -4 i a \omega ^3-b k^2+3 b \omega ^2+i (-2 c-2) \omega \right)(-4 a k \omega -2 i b k)\nonumber\\
p_3&=\left(2ik a \omega ^2-2k b \omega +2ikc\right).\label{eq183}
\end{align}
Let us deal with the factor by the term $\partial_{t_1t_1}A_0$ in (\ref{eq182}) separately.
\begin{align}
&\frac{\omega'(k)}{2ik a \omega ^2-2k b \omega +2ikc} \left(\frac{p_2}{p_3}+p_1\right)=\frac{\omega'(k)}{2ik}\frac{1}{(a\omega^2+ib\omega+c)}\left(p_1+\frac{p_2}{2ik(a\omega^2+ib\omega+c)}\right)\nonumber\\
&=\frac{\omega'(k)}{2ik}\left(\frac{p_1}{(a\omega^2+ib\omega+c)}+\frac{p_2}{2ik(a\omega^2+ib\omega+c)^2}\right)=\frac{\omega'(k)}{2ik}\left(\frac{p_1(a\omega^2+ib\omega+c)+p_2/(2ik)}{(a\omega^2+ib\omega+c)^2}\right).\label{eq184}
\end{align}
The numerator in (\ref{eq184}) becomes
\begin{align}
p_1(a\omega^2+ib\omega+c)+p_2/(2ik)&=k^2 \left(-3 a^2 \omega ^2-3 i a b \omega +a c+b^2\right)+2 a^2 \omega ^4+i a b \omega ^3+\omega ^2 (3 a-3 a c)\nonumber\\
&+\omega  (i b-2 i b c)-c^2-c.\label{eq185}
\end{align}
Now we turn once again for help to the dispersion relation (\ref{eq170}) and find
\begin{align}
a \omega ^4+i b \omega ^3+\omega ^2(c+1)&=k^2 \left(a \omega ^2+i b \omega +c\right),\nonumber\\
&\Downarrow\nonumber\\
k^2&=\frac{\omega ^2}{\left(a\omega^2+ib\omega+c\right)}+\omega^2.\label{eq186}
\end{align}
Inserting (\ref{eq186}) back into (\ref{eq185}) using $k=k$ we obtain
\begin{align}
&p_1(a\omega^2+ib\omega+c)+p_2/(2ik)=2 a^2 \omega ^4+i a b \omega ^3+\omega ^2 (3 a-3 a c)+\omega  (i b-2 i b c)-c^2-c\nonumber\\
&+\left(\frac{\omega ^2}{\left(a\omega^2+ib\omega+c\right)}+\omega^2\right) \left(-3 a^2 \omega ^2-3 i a b \omega +a c+b^2\right)\nonumber\\
&=\frac{1}{a\omega^2+ib\omega+c}\left(-a^3 \omega ^6+\omega ^4 \left(3 a b^2-3 a^2 c\right)-3 i a^2 b \omega ^5+\omega ^3 \left(-6 i a b c+i a b+i b^3\right)\right.\nonumber\\
&\left.+\omega ^2 \left(-3 a c^2+3 a c+3 b^2 c\right)-3 i b c^2 \omega -c^3-c^2\right)=\frac{p_4}{a\omega^2+ib\omega+c},\label{eq187}
\end{align}
where we denoted
\begin{align}
p_4&=-a^3 \omega ^6+\omega ^4 \left(3 a b^2-3 a^2 c\right)-3 i a^2 b \omega ^5+\omega ^3 \left(-6 i a b c+i a b+i b^3\right)\nonumber\\
&+\omega ^2 \left(-3 a c^2+3 a c+3 b^2 c\right)-3 i b c^2 \omega -c^3-c^2.\label{eq188}
\end{align}
Next we substitute (\ref{eq187}) into (\ref{eq184}) and obtain
\begin{align}
\frac{\omega'(k)}{2ik}\left(\frac{p_1(a\omega^2+ib\omega+c)+p_2/(2ik)}{(a\omega^2+ib\omega+c)^2}\right)&=\frac{\omega'(k)}{2ik}\left(\frac{p_4}{(a\omega^2+ib\omega+c)^3}\right)\nonumber\\
&=i\frac{\omega'(k)}{2k}\left(\frac{-p_4}{(a\omega^2+ib\omega+c)^3}\right).\label{eq189}
\end{align}
Using this result in (\ref{eq182}), our amplitude equation turns to
\begin{align}
\partial_{t_2}A_0+\omega'(k)\partial_{z_2}A_0+i\frac{\omega'(k)}{2k}\left(\frac{-p_4}{(a\omega^2+ib\omega+c)^3}\right)\partial_{t_1t_1}A_0-i\frac{\omega'(k)}{2k}\partial_{z_1z_1}A_0=0.\label{eq190}
\end{align}
By removing the secular terms from the equation (\ref{eq177}), the order $\varepsilon^2$ equation becomes
\begin{align}
&a^2\partial_{t_0}^{(6)}e_2+2a'b\partial_{t_0}^{(5)}e_2+(b^2-2a'-2a'c)\partial_{t_0}^{(4)}e_2-(2b'+2b'c)\partial_{t_0}^{(3)}e_2+(1+c)^2\partial_{t_0}^{(2)}e_2\nonumber\\
&-a^2\partial_{t_0}^{(4)}\partial_{z_0z_0}e_2 -2a'b\partial_{t_0}^{(3)}\partial_{z_0z_0}e_2+(2a'c-b^2)\partial_{t_0}^{(2)}\partial_{z_0z_0}e_2+2b'c\partial_{t_0}^{(1)}\partial_{z_0z_0}e_2-c^2\partial_{z_0z_0}e_2=\nonumber\\
&+NST+(*),\label{eq191}
\end{align}
where $NST$ is defined from (\ref{eq178}) as
\begin{align}
NST&=A_0^3e^{3i\theta_0}q_1+A_0^2A_0^*e^{i\theta_0}e^{2t_0\omega_i}q_2,\label{eq192}
\end{align}
and where
\begin{align}
q_1&=\left(81 a \omega ^4+27 i b \omega ^3+9 c \omega ^2\right)\nonumber\\
q_2&=3\left(a \omega ^4+8 i a \omega ^3 \omega_i-24 a \omega ^2 \omega_i^2-32 i a \omega  \omega_i^3+16 a \omega_i^4+i b \omega ^3-6 b \omega ^2 \omega_i\right.\nonumber\\
&\left.-12 i b \omega  \omega_i^2+8 b \omega_i^3+c \omega ^2+4 i c \omega  \omega_i-4 c \omega_i^2\right),\label{eq193}
\end{align}
which we solve for $e_2$ taking only the particular solution:
\begin{align}
e_2(z_0,t_0,\ldots)=c_1A_0^3e^{i3\theta_0}+c_2|A_0|^2A_0e^{i\theta_0}e^{2\omega_i}+(*),\label{eq194}
\end{align}
where
\begin{align}
c_1&=q_1/\left(\omega ^2 \left(81 a k^2-9 c-9\right)-81 a \omega ^4+27 i b k^2 \omega -27 i b \omega ^3+9 c k^2\right)\nonumber\\
&=\frac{ 9\omega^2\left(9a'\omega^2+3ib'\omega+c\right)}{-9\omega^2\left(9a'\omega^2+3ib'\omega+c\right)+9k^2\left(9a'\omega^2+3ib'\omega+c\right)-9\omega^2}\nonumber\\
&=\frac{1}{-1-\frac{1}{9a'\omega^2+3ib'\omega+c}+k^2/\omega^2}=\frac{1}{n^2(\omega)-n^2(3\omega)},\label{eq195}
\end{align}
\begin{align}
c_2&=q_2/\left(\omega ^2 \left(a k^2+24 a \omega_i^2+6 b \omega_i-c-1\right)+\omega  \left(4 i a k^2 \omega_i+32 i a \omega_i^3+i b k^2\right.\right.\nonumber\\
&\left.\left.+12 i b \omega_i^2-4 i c \omega_i-4 i \omega_i\right)+\omega ^3 (-8 i a \omega_i-i b)-4 a k^2 \omega_i^2-a \omega ^4\right.\nonumber\\
&\left.-16 a\omega_i^4-2 b k^2 \omega_i-8 b \omega_i^3+c k^2+4 c \omega_i^2+4 \omega_i^2\right)\nonumber\\
&=\frac{3\left(\omega^*\right)^2\left(a\left(\omega^*\right)^2+ib\omega^*+c\right)}{k^2 \left(a \left(\omega^*\right)^2+i b \omega^*+c\right)-a \left(\omega^*\right)^4-i b \left(\omega^*\right)^3+(-c-1) \left(\omega^*\right)^2}\nonumber\\
&=\frac{3\left(\omega^*\right)^2\left(a\left(\omega^*\right)^2+ib\omega^*+c\right)}{k^2 \left(a \left(\omega^*\right)^2+i b \omega^*+c\right)-\left(\omega^*\right)^2\left( a \left(\omega^*\right)^2+i b \left(\omega^*\right)+c\right)- \left(\omega^*\right)^2}\nonumber\\
&=\frac{3\left(\omega^*\right)^2}{k^2-\left(\omega^*\right)^2- \left(\omega^*\right)^2/\left(a\left(\omega^*\right)^2+ib\omega^*+c\right)}=\frac{3\left(\omega^*\right)^2}{k^2-\left(\omega^*\right)^2\left(1+ 1/\left(a\left(\omega^*\right)^2+ib\omega^*+c\right)\right)}\nonumber\\
&=\frac{3\left(\omega^*\right)^2}{k^2-\left(\omega^*\right)^2n^2\left(\omega^*\right)},\label{eq196}
\end{align}
where $\omega^*=\omega+2i\omega_i$.

Defining as before the amplitude as in (\ref{eq45}) and proceeding the usual way using (\ref{eq174}) and (\ref{eq190}) we get the final amplitude equation
\begin{align}
\partial_tA+\omega'(k)\partial_zA-i\beta\partial_{zz}A+i\alpha\partial_{tt}A=0,\label{eq197}
\end{align}
where
\begin{align}
\alpha&=\frac{\omega'(k)}{2k}\left(\frac{-p_4}{(a\omega^2+ib\omega+c)^3}\right),\label{eq198}\\
\beta&=\frac{\omega'(k)}{2k},\label{eq199}
\end{align}
and where $p_4$ is defined in (\ref{eq188}). The overall approximate solution to (\ref{eq164}) is then
\begin{align}
E(z,t)&=A(z,t)e^{i(kz-\omega t)}+c_1\varepsilon^2A^3(z,t)e^{i3(kz-\omega t)}\nonumber\\
&+c_2\varepsilon^2|A(z,t)|^2A(z,t)e^{i(kz-\omega t)}e^{2t\omega_i}+(*),\label{eq200}
\end{align}
where $c_1,c_2$ are defined in (\ref{eq195}) and (\ref{eq196}). To verify that $\alpha,\beta$ in (\ref{eq39}) and (\ref{eq198}) are the same, we look at (\ref{eq39}) and get
\begin{align}
\alpha&=\omega'(k)\frac{n^2(\omega)+2\omega \sqrt{2\pi}\hat{\chi}'(\omega)+\frac{1}{2}\omega^2\sqrt{2\pi}\hat{\chi}''(\omega)}{2k}\nonumber\\
&=\frac{\omega'(k)}{2k}\left(1+\frac{1}{a\omega^2+ib\omega+c}-2\omega\frac{2 a \omega +i b}{\left(a \omega ^2+i b \omega +c\right)^2}-\frac{1}{2}\omega^2 \frac{2 \left(-3 i a b \omega +a \left(c-3 a \omega ^2\right)+b^2\right)}{(a\omega^2+ib\omega+c)^3}\right)\nonumber\\
&=\frac{\omega'(k)}{2k}\frac{1}{(a\omega^2+ib\omega+c)^3}\left(a^3 \omega ^6+\omega ^4 \left(3 a^2 c-3 a b^2\right)+3 i a^2 b \omega ^5+\omega ^3 \left(6 i a b c-i a b-i b^3\right)\right.\nonumber\\
&\left.+\omega ^2 \left(3 a c^2-3 a c-3 b^2 c\right)+3 i b c^2 \omega +c^3+c^2\right).\label{eq201}
\end{align}
We can see that the expression in the numerator is exactly the same as (\ref{eq188}) multiplied with -1.

Using the same argument as in (\ref{eq48}) we can simplify the amplitude equation (\ref{eq197}) as
\begin{align}
\partial_tA+\omega'(k)\partial_zA-i\partial_{zz}A\left(\beta-\alpha\left(\omega'(k)\right)^2\right)=0,\label{eq202}
\end{align}
which can be solved as an initial value problem.

\end{appendices}

\bibliographystyle{unsrt}
\bibliography{paper}

\end{document}